\documentclass[11pt,twoside,a4paper]{article}

\usepackage{amsmath,bbm,amssymb,mathtools,geometry, xcolor,mathabx,mathrsfs,amsthm}
\allowdisplaybreaks
\usepackage[margin=1in]{caption}

\usepackage[usenames,dvipsnames]{pstricks}
 \usepackage{pstricks-add}
 \usepackage{epsfig}

\mathtoolsset{showonlyrefs}

\newcommand{\D}{\nabla}
\newcommand{\dd}{\partial}
\newcommand{\norm}[1]{\left\lVert#1\right\rVert}

\newtheorem{Thm}{Theorem}[section]
\newtheorem{Lem}[Thm]{Lemma}
\newtheorem{Rem}[Thm]{Remark}

\newtheorem{Def}[Thm]{Definition}
\newtheorem{Prop}[Thm]{Proposition}

\newtheorem{Cor}[Thm]{Corollary}

\def\R{\mathbb R}

\def\N{\mathbb N}

\newcommand{\mres}{\mathbin{\vrule height 1.6ex depth 0pt width
0.13ex\vrule height 0.13ex depth 0pt width 1.3ex}}

\def\Xint#1{\mathchoice
   {\XXint\displaystyle\textstyle{#1}}
   {\XXint\textstyle\scriptstyle{#1}}
   {\XXint\scriptstyle\scriptscriptstyle{#1}}
   {\XXint\scriptscriptstyle\scriptscriptstyle{#1}}
   \!\int}
\def\XXint#1#2#3{{\setbox0=\hbox{$#1{#2#3}{\int}$}
     \vcenter{\hbox{$#2#3$}}\kern-.5\wd0}}

\def\dashint{\Xint-}
\newcommand{\prodscal}[2]{\langle{#1},{#2}\rangle}

\usepackage{appendix}

\title{Quantitative regularity properties for the optimal design problem}
\author{Lorenzo Lamberti\footnote{lorenzo.lamberti@univ-lorraine.fr, Université de Lorraine, CNRS, IECL, F-54000 Nancy, France} \; and Antoine Lemenant\footnote{antoine.lemenant@univ-lorraine.fr, Université de Lorraine, CNRS, IECL, F-54000 Nancy, France}}

\begin{document}

\maketitle

\abstract{In this paper we slightly improve the regularity theory for the so called optimal design problem. We first establish the uniform rectifiability of the  boundary of the optimal set, for a larger class of minimizers, in any dimension. As an application, we improve the bound obtained by Larsen in dimension~2 about the mutual distance between two connected components. Finally we also prove that the full regularity in dimension 2 holds true provided that the ratio between the two constants in front of the Dirichlet energy is not larger than 4, which partially answers to a question raised by Larsen.  }

\tableofcontents

\section{Introduction}

Let $\Omega\subset\R^N$ be a bounded connected open set and $0<\alpha<\beta$ be two constants. We define $\sigma_E=\alpha\mathbbm{1}_E+\beta\mathbbm{1}_{\Omega\setminus E}$. The so called \emph{optimal design problem} consists of minimizing among couples $(u,E)$ the following problem
\begin{equation}
\min_{(u,E) \in \mathcal{A}}\int_{\Omega} \sigma_{E} |\nabla u|^2 \;dx + P(E; \Omega), \label{optimalDesign}
\end{equation}
where
$$\mathcal{A}:=\left\{(u,E) \text{ s.t. } |E|=V_0 \text{ and }u\in H^1(\Omega), u=u_0 \text{ on } \partial \Omega\right\}.$$
Here  $V_0\in(0,|\Omega|)$ is a given volume and    $u_0 \in H^1(\Omega)$ is  a boundary datum in the sense that $u=u_0$ on $\partial \Omega$ means $u-u_0 \in H^1_0(\Omega)$.

This problem has been widely studied by many famous authors from the 90's up to nowadays (see for instance \cite{AmbrosioButtazzo,  DePhilippisFigalli, espFusco, EspLam2,  FuscoJulin, Larsen1999,Larsen2003,  Lin}), and a lot is known about the regularity of minimizers.

To provide some historical context, in 1993 Ambrosio and Buttazzo \cite{AmbrosioButtazzo} established the existence of solutions together with the higher integrabilty of the gradient of the deformation $u$. In the same year, Lin \cite{Lin} proved that for a minimizer $(u,E)$, the function $u$ must be globally $C^{0,1/2}$-regular in $\Omega$ and the boundary of $E$ inside $\Omega$ is $C^{1,\alpha}$ outside a singular set of  zero $\mathcal{H}^{N-1}$-measure.  This was again proved with different techniques by Fusco and Julin \cite{FuscoJulin} in 2015, and improved in the sense that the singular set must have Hausdorff dimension strictly less than $N-1$. In the same year, De Philippis and Figalli \cite{DePhilippisFigalli} independently obtained the estimate about the Hausdorff dimension of the singular set, by employing porosity techniques.  
Later, some variants of the problems with more general densities with quadratic growth and $p$-growth, higher-order operators or in a vectorial context have been studied in \cite{Arroyo-Rabasa, CarEspLam1, CarEspLam2, CFP, CFP2, Esp, EspLam, EspLamPis, Lam, LinKohn}.

However, in the specific dimension $N=2$, the full regularity of $E$ is still a challenging open problem, raised by Larsen. Indeed, Larsen \cite{Larsen2003} proved that any connected component of $E$ has a $C^1$ boundary. This does not prevent that a countable number of connected components accumulate in a way that $\partial E$ may not be globally smooth, but Larsen conjectures (in \cite{Larsen1999} and again in \cite{Larsen2003}) that it might not be the case.

The first result of this paper is a positive answer to Larsen's conjecture in the case when $\beta< 4 \alpha$. Here is our first result.

\begin{Thm}\label{main1} Let $N=2$ and  let $(u,E)$ be a minimizer for the optimal design Problem \eqref{optimalDesign}. Assume moreoever  that $\beta< 4 \alpha$.  Then $\partial E$ is a smooth $C^{1,\alpha}$-surface in $\Omega$. 
\end{Thm}

The proof of Theorem  \ref{main1} is very short and is given in Section \ref{sectionsmooth}. It relies on a monotonicity formula, similar to that of Bonnet \cite{bonnet}, which directly establishes that $E$ is an almost minimizer for the perimeter in the regime $\beta< 4 \alpha$, allowing us to apply the standard regularity theory. In the same section, we investigate some further monotonicity properties that imply, for instance, that, without any restriction on $\alpha,\beta$, if $\partial E$ intersects $\partial B_s(x)$ by only two points for all $s\in(0,r)$, then $\partial E\cap B_{\frac{r}{2}}(x)$ must be smooth.

Theorem \ref{main1} also partially improves an earlier result of Esposito and Fusco \cite{espFusco}, in which they prove that $\dd E$ is a smooth surface when $\beta\leq \gamma_N \alpha$, where $\gamma_N$ is an explicit constant depending on dimension. In particular, for $N=2$, they obtain the constant $\gamma_2=5/3$. Since this value is strictly less than $4$,  Theorem \ref{main1} is an improvement of \cite[Theorem 2]{espFusco}, in the special case of $N=2$.

Notice that a possible way to solve Larsen's conjecture would be to prove that $E$ admits a finite number of connected components.  Subsequently, any qualitative information about the connected components of $E$ would be of great interest.  Toward this direction Larsen was able to prove in \cite{Larsen1999} that for two given connected components $E_1,E_2$ of $E$, it holds that
\begin{equation}
{\rm dist}(E_1,E_2)>0, \label{distanceEq}
\end{equation}
which  actually stands for the main result of \cite{Larsen1999}.

Our second main result is a quantitative improvement of Larsen's estimate \eqref{distanceEq}. In the following statements, there is no more restrictions on the values of  $\alpha$ and $\beta$.

\begin{Thm}\label{main2}  Let $N=2$ and  let $(u,E)$ be a minimizer for the optimal design Problem \eqref{optimalDesign}.  Then there exist two constants $C_0>0$ and $\varepsilon_0>0$ such that  for any two components $E_1,E_2$ of $E$ it holds that
 \begin{equation*}
    \text{ either }  \quad \mathrm{dist}(E_1,E_2)\geq \varepsilon_0 \quad\text{or}\quad {\rm dist}(E_1,E_2)^2\geq C_0 \min\{|E_1|,|E_2|\}.
 \end{equation*}
\end{Thm}

The proof of Theorem \ref{main2} is given in Section \ref{SectionDist}, and uses the   uniform rectifiability of $\partial E$. This fact is   established first  in Section \ref{uniformRect}, in a much more general context, and it is interesting for its own.

Indeed, a wide proportion of the present paper is to prove the uniform rectifiability of $\partial E$, which is actually valid in any dimension, and applies to the more general class of quasi-minimizers in the sense of David and Semmes. We also relax the volume constraint by working  with a penalized version of the functional. As it was shown by Esposito and Fusco \cite{espFusco}, the minimization problem with this penalized functional is equivalent to the original problem (in the language of \cite{espFusco}, we work with a generalization of $\Lambda$-minimizers).   More precisely, we introduce 
$$\mathcal{A'}:=\left\{(u,E) \text{ s.t. }  E\subset \Omega \text{ and }u\in H^1(\Omega), u=u_0 \text{ on } \partial \Omega\right\},$$
and then consider the following problem:

\begin{equation}
\min_{(u,E) \in \mathcal{A'}} \mathcal{F}(u,E),\label{optimalDesign2}
\end{equation}
with 
\begin{equation}
 \mathcal{F}(u,E) :=\int_{\Omega} \sigma_{E} |\nabla u|^2 \;dx + \Psi_E(\Omega) + \Lambda||A|-V_0|.\label{optimalDesign22}
\end{equation}
Here $\Lambda>0$  is a constant, and $\Psi_E$ is a scalar Radon measure that we assume to be comparable to the perimeter, that is
\begin{equation}
\label{Comparability}
K^{-1}P(E; G)\leq \Psi_E(G)\leq K P(E; G),
\end{equation}
for any set $G\subset\R^N$ and for some constant $K>0$. In the case when $K=1$, we recover the so-called $\Lambda$-minimizers of the classical optimal design problem. It was furthermore proved in \cite[Theorem 1]{espFusco} that minimizers of the constrained Problem \eqref{optimalDesign} are also $\Lambda$-minimizers for a suitable choice of $\Lambda>0$ (see also Theorem \ref{EquivPen}).

Here is the regularity result that we obtain with regards to $\Lambda$-quasi-minimizers (that in the sequel will be sometimes simply called \emph{quasi minimizers}).


\begin{Thm}\label{main3} Let $(u,E)$ be a $\Lambda$-quasi minimizer for the optimal design problem, i.e. a minimizer for the Problem \eqref{optimalDesign2}. Then $E$ satisfies the condition-B (see Definition \ref{defconditionB}). In particular, $\partial E$ is uniformly rectifiable in $\Omega$. 
\end{Thm}

The notion of uniform rectifiabilty is a sort of quantitative notion of rectifiability that was introduced and intensively studied by David and Semmes (see for instance \cite{david-semmes}). In particular, it provides some nice uniform control in all scales, such as big pieces of Lipschitz graphs, smallness of the flatness in many balls in a uniform way, etc. This notion is more global and quantitative on the whole set $\partial E$ compared to the usual standard local regularity results such as $\varepsilon$-regularity type ones. The combination of uniform rectifiability and local regularity gives rise to new interesting statements.

As already pointed out, Theorem \ref{main2} is an example of those statements that use the uniform rectifiability of $\partial E$. In addition, we get several other consequences of uniform rectifiability, such as a new way to improve the Hausdorff dimension of the singular set (see Corollary \ref{hausdorff}) different from \cite{DePhilippisFigalli} and  \cite{FuscoJulin}. A last example is  Proposition \ref{Larsen4} that gives an estimate already obtained before by Larsen in \cite{Larsen2003}, for which we provide here  a completely different proof relying on the uniform rectifiability of $\partial E$.

To prove Theorem \ref{main3} we first show that quasi-minimizers are Ahlfors-regular, adapting the standard proof already known for optimal design minimizers. Then we prove that $E$ satisfies the so-called ``condition-B'' (see Definition \ref{defconditionB}). For that purpose we use a control of the normalized energy of $u$ by Carleson measure estimates. The uniform rectifiability follows immediately, as it is known from David and Semmes \cite{david-semmes} that it is a consequence of Ahlfors-regularity and condition-B.

Once the uniform rectifiability is established, we use it in dimension $N=2$ to prove Theorem \ref{main2}. The general strategy follows the original one of Larsen \cite{Larsen1999} in his proof of \eqref{distanceEq}, but incorporating the uniform rectifiability to make it quantitative. Also, we take benefit from this paper to entirely rewrite  the original arguments of Larsen, especially the key Lemma \ref{LemmaRectangle} for which we used some of his ideas, but written here with completely different arguments that we believe are more detailed than what can be found in \cite{Larsen1999}.\\





{\bf Acknowledgement.} This paper was partially financed by the junior IUF grant of A. Lemenant and by the ANR project ``STOIQUES''. Lorenzo Lamberti is a member of the Gruppo Nazionale per l’Analisi Matematica, la Probabilità e le loro Applicazioni (GNAMPA) of the Istituto Nazionale di Alta Matematica (INdAM).

\section{Notation and preliminary definitions}
Let $\Omega$ be a bounded connected open subset of $\R^N$, with  $N\geq 2$. We denote by 
$
B_r(x):=\left\{y\in \R^N: |y-x|<r\right\}$  the open ball centered at $x\in \R^N$ of radius $r>0$ and  as usual $\omega_N$ stands for the Lebesgue measure of the unit ball in $\R^N$. If $x_0=0$ we simply write $B_r$. We denote by $\sharp E$ the cardinality of the set $E$ and by $C$ a
generic constant that may vary from line to line. We write $\prodscal{\xi}{\eta}$ for the inner product of vectors $\xi, \eta \in \mathbb{R}^N$, and consequently $|\xi|:=\prodscal{\xi}{\xi}^{\frac{1}{2}}$ will be the corresponding Euclidean norm. In the following, we  denote
\begin{equation*}
\mathbf{C}_{r}(x_0):=x_0+\{y\in\R^N\,:\,|y_N|<r,\,|y-y_Ne_N|<r\},
\end{equation*}
the cylinder centered in $x_0\in\R^N$ with radius $r>0$ oriented in the direction of the $N$-th versor $e_N$. For $i\in\{1,\dots,N\}$, we call by $\pi_i$ the projection on the $i$-th coordinate, i.e. $\pi_i(x)=x_i$, for $x=(x_1,\dots,x_N)\in\R^N$.\\
\indent Let $G\subset\R^N$. We define the set of points of density $t\in[0,1]$ as follows:
\begin{equation*}
    G^{(t)}=\left\{x\in \R^N\,:\, \lim_{r\rightarrow 0^+}\frac{|G\cap B_r(x)|}{|B_r(x)|}=t\right\}.
\end{equation*}
Let $U$ be an open subset of $\R^n$. A Lebesgue measurable set $E\subset \R^N$
is said to be a set of locally finite perimeter in $U$ if there exists a $\R^N$-valued Radon measure $\mu_E$ on $U$ (called the Gauss-Green measure of $E$) such that
\begin{equation*}
\int_{E}\nabla \phi\ dx=\int_{U}\phi \, d\mu_E,\quad\forall \phi \in C^1_c(U).
\end{equation*}
Moreover, we denote the perimeter of $E$ relative to $G\subset U$ by $P(E,G)=|\mu_E|(G)$.\\
 \noindent It is well known that the support of $\mu_E$ can be characterized by
\begin{equation}\label{support}
\text{spt}\mu_E=\bigl\{x\in U: 0<|E\cap B_r(x)|<\omega_n r^n, \,\forall r>0\bigr\}\subset U\cap \partial E,
\end{equation}
(see \cite[Proposition 12.19]{maggi}).
If $E$ is of finite perimeter in $U$, the {\it reduced boundary} $\partial^*E \subset U$ of $E$ is the set of those $x\in U$ such that
\begin{equation}\label{RB}
\nu_E(x):=\lim_{r\rightarrow 0^+}\frac{\mu_E(B_r(x))}{|\mu_E|(B_r(x))}
\end{equation}
exists and belongs to $\mathbb S^{n-1}$. We address the reader to \cite{maggi} for a complete dissertation about sets of finite perimeter.\\
\indent For $u\in H^1(B_r(x_0))$ and $p\in[1,2]$ we denote
\begin{equation*}
    \omega_p(x_0,r)=r^{1-\frac{2N}{p}}\bigg(\int_{B_r(x_0)}|\D u|^p\,dx\bigg)^{\frac{2}{p}}.
\end{equation*}
We simply write $\omega(x_0,r):=\omega_2(x_0,r)$. In the subsequent sections, we need the definition of Alhfors-regular sets.

\begin{Def}[Alhfors-regularity] Let $G\subset\R^N$ be a closed set. We say that $G$ is $(N-1)$-Ahlfors-regular (or, shortly, Ahlfors-regular) if there exists a positive constant $C_A$ such that
\begin{equation*}
    C_A^{-1} r^{N-1}\leq\mathcal{H}^{N-1}(G\cap B_r(x_0))\leq C_A r^{N-1},\quad\forall x_0\in G,\,\forall r>0.
\end{equation*}
   
\end{Def}

In what follows, the definition of uniformly rectifiable set will be needed. It is a stronger and more quantitative notion of rectifiability. There are many equivalent (and not simple) definitions of uniform rectifiability.  For instance, here is one of them.

\begin{Def}[Uniform rectifiability] \label{Uniform_rectifiability}
Let $G\subset\R^N$ be an Ahlfors-regular set. We say that $G$ is uniformly rectifiable if there exist two positive constants $\theta$ and $C$ such that, for each ball $B$ centered in $G$, we can find a compact set $A\subset\R^{N-1}$ and a bi-Lipschitz map $\rho\colon A\rightarrow \R^N$ such that
\begin{equation*}
    C^{-1}|x-y|\leq |\rho(x)-\rho(y)|\leq C|x-y|,\quad\forall x,y\in A,
\end{equation*}
and
\begin{equation}
    \mathcal{H}^{N-1}(G\cap\rho(A)\cap B)\geq \theta\mathcal{H}^{N-1}(G\cap B).
\end{equation}
\end{Def}

Uniformly rectifiable sets have been extensively studied in the monography \cite{david-semmes}. For example, they provide a connection between geometric measure theory and harmonic analysis.  In this paper, we shall make use of a geometric characterization of uniform rectifiability. We first need the following definition.

\begin{Def}[Carleson sets]
    Let $G\subset\R^N$ be Alhfors-regular. We say that a measurable set $A\subset G\times \R_+$ is a Carleson set if $\mathbbm{1}_Ad\mathcal{H}^{N-1}\mres G\frac{dt}{t}$ is a Carleson measure on $G\times\R_+$, i.e. there exists a positive constant $C$ such that
    \begin{equation*}
        \int_0^r\int_{G\cap B_r(z)}\mathbbm{1}_A(x,t)\,d\mathcal{H}^{N-1}\frac{dt}{t}\leq C r^{N-1},\quad\forall z\in G,\,r>0.
    \end{equation*}
\end{Def}
This is an invariant way of saying that the set $A$ is enough small and that it behaves as it was $(N-1)$-dimensional from the perspective of $G\times\{0\}$.

Here it follows a useful characterization of uniform rectifiabilty that can be found in \cite[Theorem 2.4]{david-semmes}.
\begin{Prop}
\label{EquivBWGL}
    Let $G\subset\R^N$ be an Ahlfors-regular set. Then, $G$ is uniformly rectifiable if and only if $G$ satisfies the bilateral weak geometric lemma (BWGL), i.e., for every $\varepsilon>0$,
    \begin{equation*}
        \{(x,t)\in G\times\R_+\,:\, \beta(x,t)>\varepsilon\}
    \end{equation*}
    is a Carleson set. Here, the quantity
    \begin{equation*}
\beta(x,t):=\inf_{\substack{P\subset\R^N \\ P \text{affine hyperplane}}}\bigg\{\sup_{y\in G\cap B_t(x)}t^{-1}\mathrm{dist}(y,P)+\sup_{z\in P\cap B_t(x)}t^{-1}\mathrm{dist}(z,G)\bigg\} 
    \end{equation*}
    denotes the bilateral flatness at the point $x\in G$ at scale $t\in\R_+$.
\end{Prop}

This equivalence allows us to have a quantitative control of the flatness. In other words, the sets of points where the flatness of the set $E$ is arbitrarily small is big in terms of measure. This ensures the existence of many balls centered in the boundary of the optimal shape where the well-known result of $\varepsilon$-regularity holds.\\

\indent In practice, it is not so easy to prove uniform rectifiability from Definition \ref{Uniform_rectifiability} or Proposition~\ref{EquivBWGL}. For the particular case  of boundaries of sets, there exists a nice criterium using the so-called condition-B, which we present in the next definition.


\begin{Def}[condition-B]\label{defconditionB}
Let $G$ be a measurable subset of $\R^N$. We say that $G$ satisfies the condition $B$ in $\Omega$ if $G$ is open, $\dd G$ is Ahlfors-regular and if for any open set $U\subset\subset\Omega$ of $\R^N$ there exist two constants $C_0>1$ and $r_0\in(0,\mathrm{dist}(U,\dd\Omega))$ such that for any $x_0\in\dd G$ and $r\in (0,r_0)$, we can find two balls $B_1\subset B_r(x_0)\cap E$ and $B_2\subset B_r(x_0)\setminus\overline{E}$ with radius greater than $C_0^{-1}r$.
\end{Def}


The next proposition follows by combining \cite[Theorem 1.20, Proposition 1.18, Theorem 1.14 and Proposition 3.35]{david-semmes2}.

\begin{Prop}
\label{CritUnifRect}
    Let $G\subset\R^N$ be an open such that $\partial G$ is an Ahlfors-regular set. If furthermore $G$  satisfies the condition-B, then $\partial G$ is uniformly rectifiable.
\end{Prop}

To conclude the section, we cite the following theorem, whose proof is contained in \cite{EspLamPis}.

\begin{Thm}[\cite{EspLamPis}]
\label{EquivPen}
    There exists a constant $\Lambda_0>0$ such that if $(u,E)$ is a minimizer of the functional
    \begin{equation}
\label{FunctPen}\int_\Omega\sigma_F|\D w|^2\,dx+\Psi_F(\Omega)+\Lambda||F|-V_0|
    \end{equation}
    for some $\Lambda\geq\Lambda_0$ among all the configurations $(F,w)$ such that $w=u_0$ on $\dd\Omega$, then $|E|=V_0$ and $(E,u)$ is a minimizer of Problem \eqref{optimalDesign2}. Conversely, if $(E,u)$ is a minimizer of Problem \eqref{optimalDesign2}, then it is a minimizer of \eqref{FunctPen}, for any $\Lambda>0$.
\end{Thm}

\section{Uniform rectifiabilty for quasi-minimizers in dimension~$N$}
\label{uniformRect}

This section is devoted to prove that the boundary of the optimal set is uniformly rectifiable. In view of this aim, we first show that it is Alhfors-regular. Afterwards, it suffices to prove that it satisfies the condition-B (see Proposition \ref{CritUnifRect}). We show the validity of the latter property in Proposition \ref{rigot}. \\
\indent Throughout the entire section we assume that $E$ is a Borel set with

\begin{equation*}
\dd E=\mathrm{spt}(\mu_E)=\{x\in\R^N\,:\, 0<|E\cap B_r(x)|<|B_r(x)|, \,\forall r>0\}.
\end{equation*}

Let us emphasize that, according to \cite[Proposition 12.19]{maggi}, for any open set of finite perimeter, one can always find an equivalent Borel set with this property. Furthermore, it is easy to show that $E^{(1)}$ is a valid choice. At the end of this section, we prove that for a quasi-minimizer, one can actually choose this Borel set to be an open set (see Lemma \ref{rigot}).\\

\indent In the following theorem, we prove that the boundary of a minimal set is Alhfors-regular. The scheme of proof is rather standard: it follows the original proof of Ahlfors-regularity for the optimal design problem that we adapt for quasi-minimizers instead of minimizers. We rewrite shortly the proof for the convenience of the reader.

\begin{Thm}[Ahlfors regularity]
Let $(u,E)$ be a minimizer of \eqref{optimalDesign2} and $U\subset\subset\Omega$ be an open set. Then there exists a positive constant $C_A=C_A(N,\alpha,\beta,\Lambda,K,\norm{\D u}_{L^2(\Omega)}\big)$ such that, for every $x_0\in\dd E$ and $B_r(x_0)\subset U$, it holds that
\begin{equation}
\label{DensityPerimeter}
\frac{1}{C_A}r^{N-1}\leq P(E;B_r(x_0))\leq C_A r^{N-1}.
\end{equation}
Furthermore, $\mathcal{H}^{N-1}((\dd E\setminus \dd^*E)\cap\Omega)=0$ and $\dd E$ is Ahlfors-regular.
\end{Thm}

\begin{proof}The proof is divided in four steps.\\ \indent \textbf{Step 1:} \emph{Upper bound on the energy.} We show that for every open set $U\subset\subset \Omega$ there exists a constant $C=C(N,\alpha,\beta,K,\Lambda)>0$ such that for every $B_r(x_0)\subset U$ it holds
\begin{equation}
\label{eqq11}
\mathcal{F}(E,u;B_r(x_0))\leq C r^{N-1}.
\end{equation}
In order to prove it, using the $\Lambda$-minimality of $(u,E)$ with respect to $(u,E\cup B_r(x_0))$ (see Theorem \ref{EquivPen}) and the comparability condition \eqref{Comparability}, one can obtain
\begin{equation}
\label{eqq19}
(\beta-\alpha)\int_{B_r(x_0)\setminus E}|\D u|^2\,dx + \frac{1}{K}P(E;B_r(x_0))\leq C r^{N-1},
\end{equation}
where $C=C(N,K,\Lambda)$.
To show \eqref{eqq11} it suffices to prove that there exist some constants $M>0$,  $\tau\in\big(0,\frac{1}{2}\big)$ and $h_0\in\N$, depending on $N$, $\frac{\beta}{\alpha}$ and $K$,  such that, for any $B_r(x_0)\subset U$, we have
\begin{equation*}
\int_{B_r(x_0)}|\D u|^2\leq h_0r^{N-1} \quad\text{or}\quad \int_{B_{\tau r}(x_0)}|\D u|^2\,dx\leq M\tau^{N-\frac{1}{2}}\int_{B_r(x_0)}|\D u|^2\,dx.
\end{equation*}
We assume by contradiction that for some $M>0$ and $\tau\in\big(0,\frac{1}{2}\big)$ to be chosen, for any $h\in\N$ there exists a ball $B_r(x_h)\subset U$ such that 
\begin{equation}
\label{eqq23}
\int_{B_{r_h}(x_h)}|\D u|^2\,dx\geq hr^{N-1} \quad\text{and}\quad \int_{B_{\tau r_h}(x_h)}|\D u|^2\,dx\geq M\tau^{N-\frac{1}{2}}\int_{B_{r_h}(x_h)}|\D u|^2\,dx.
\end{equation}
Combining the first inequality 
 and \eqref{eqq19}, we get
\begin{equation}
\label{a41}
\int_{B_{r_h}(x_h)\cap E}|\D u|^2\,dx<\frac{C}{h}\int_{B_{r_h}(x_h)}|\D u|^2\,dx,
\end{equation}
where $C=C\big(N,\frac{\alpha}{\beta},K,\Lambda)$.
For $y\in B_1$, we define
\begin{equation*}
v_h(y):=\frac{u(x_h+r_hy)-a_h}{\varsigma_h r_h},
\end{equation*}
where we have denoted 
\begin{equation*}
a_h:=\dashint_{B_{r_h}(x_h)}u\,dx \quad\text{and}\quad\varsigma_h^2:=\dashint_{B_{r_h}(x_h)}|\D u|^2\,dx.
\end{equation*}
Furthermore we set
\begin{equation*}
E^*_h:=B_1\setminus \frac{E-x_h}{r_h}.
\end{equation*}
Since $\{\D v_h\}_{h\in\N}$ is bounded in $L^2(B_1)$, there exist a (not relabeled) subsequence of $v_h$ and $v\in H^1(B_1)$ such that $v_h\rightharpoonup v$ in $H^1(B_1)$ and $v_h\rightarrow v$ in $L^2(B_1)$. Furthermore, using the upper bound on the perimeters of $E^*_h$ in $B_1$ given by \eqref{eqq19}, up to a not relabeled subsequence, $\mathbbm{1}_{E^*_h}\rightarrow\mathbbm{1}_{E^*}$ in $L^1(B_1)$, for some set $E^*\subset B_1$ of locally finite perimeter.\\
\indent From the minimality of $u$, we obtain the following minimality relation for $v_h$:
\begin{equation}
\label{eqq22}
\int_{B_1}\sigma_{B_1\setminus E^*_h}|\D v_h|^2\,dx\leq \int_{B_1}\sigma_{B_1\setminus E^*_h}\big|\D v_h+\varsigma_h^{-1}\D\psi\big|^2\,dy,\quad\forall\psi\in H^1(B_1).
\end{equation}
Choosing $\psi_h=\varsigma_h \eta (v-v_h)$, where $\eta\in C^1_c(B_1)$, with $0\leq\eta\leq 1$, we get
\begin{align*}
&\int_{B_1}\sigma_{B_1\setminus E_h}|\D v_h|^2\,dy\leq 
\int_{B_1}\sigma_{B_1\setminus E^*_h}|\eta\D v+(1-\eta)\D v_h|^2\,dy\\
& + \int_{B_1}\sigma_{B_1\setminus E^*_h}(v-v_h)^2|\D\eta|^2\,dy+2\int_{B_1}(v-v_h)\prodscal{\D\eta}{\eta\D v+(1-\eta)\D v_h}\,dy\\
& \leq \int_{B_1}\sigma_{B_1\setminus E^*_h}\eta|\D v|^2\,dy+\int_{B_1}\sigma_{B_1\setminus E^*_h}(1-\eta)|\D v_h|^2\,dy+o(1),
\end{align*}
where we have used the convergence of $v_h$ and the boundedness of $\{\D v_h\}_{h\in\N}$. Thus, we obtain
\begin{align}
\label{eqq21}
\int_{B_1}\sigma_{B_1\setminus E^*_h}\eta|\D v_h|^2\,dy\leq \int_{B_1}\sigma_{B_1\setminus E^*_h}\eta|\D v|^2\,dy+o(1).
\end{align}
Furthermore, by \eqref{a41} and the and the equi-integrability of $\{\D v_h\}_{h\in\N}$, we deduce that
\begin{equation}
\label{eqq20}
    \lim_{h\rightarrow +\infty}\int_{E^*_h}|\D v_h|^2\,dy=0 \quad\text{and}\quad \int_{E^*}|\D v|^2\,dy=\lim_{h\rightarrow +\infty}\int_{E^*_h}|\D v|^2\,dy=0.
\end{equation}
Thus, we may rewrite \eqref{eqq21} as follows:
\begin{align}
\int_{B_1\setminus E^*_h}\eta|\D v_h|^2\,dy\leq \int_{B_1\setminus E^*_h}\eta|\D v|^2\,dy+o(1).
\end{align}
Passing to the upper limit as $h\rightarrow+\infty$, using the lower semicontinuity and letting $\eta\rightarrow 1$, we get
\begin{align*} \lim_{h\rightarrow+\infty}\int_{B_1\setminus E^*}|\D v|^2\,dy=\int_{B_1\setminus E^*}|\D v|^2\,dy.
\end{align*}
Using also the second equality in \eqref{eqq20}, we infer that $\D v_h\rightarrow\D v$ in $L^2(B_1)$ and therefore $v_h\rightarrow v$ in $H^1(B_1)$. Letting $h\rightarrow+\infty$ in \eqref{eqq22}, we infer that $v$ minimizes
\begin{equation*}
\int_{B_1}\sigma_{B_1\setminus E^*}|\D v|^2\,dy.
\end{equation*}
Thus, there exist two constants $\tau_0\in\big(0,\frac{1}{2}\big)$ and $\overline{C}>0$ such that
\begin{equation*}
    \dashint_{B_\tau}|\D v|^2\,dy\leq \overline{C}\dashint_{B_1}|\D v|^2\,dy=\overline{C}\lim_{h\rightarrow+\infty}\dashint_{B_1}|\D v_h|^2\,dx=\overline{C}\omega_n.
\end{equation*}
In conclusion, choosing $M>\overline{C}\omega_n$, by \eqref{eqq23} we get
\begin{align*}
    \int_{B_\tau}|\D v|^2\,dy< M\leq M\tau^{-\frac{1}{2}}\leq \int_{B_\tau}|\D v|^2\,dy,
\end{align*}
which is a contradiction.\\
\indent \textbf{Step 2:} \emph{Decay of the energy in the balls where the perimeter of $E$ is small.}
We want to show that for every $\tau\in (0,1)$ there exists $\varepsilon_0=\varepsilon_0(\tau)>0$ such that, if $B_r(x_0)\subset \Omega$ and $P(E;B_r(x_0))<\varepsilon_0 r^{N-1}$, then
\begin{equation}\label{D}
\mathcal F(E,u;B_{\tau r}(x_0))\leq C \tau^N\bigl(\mathcal F(E,u;B_r(x_0))+r^N\bigr),
\end{equation}
for some positive constant $C=C\big(N,\alpha,\beta,\Lambda,K,\norm{\D u}_{L^2(\Omega)}\big)>0$ independent of $\tau$ and $r$.
First of all, we remark that under the assumption \eqref{Comparability}, $\Psi_E$ is absolutely continuous with respect to $\mathcal{H}^{N-1}\mres\dd^*E$. Therefore, by the Radon-Nikodym Theorem there exists a function $\theta\colon \Omega\rightarrow\R$ such that
\begin{equation}
\Psi_E(G)=\int_{\dd^*E \cap G}\theta\,d\mathcal{H}^{N-1},
\end{equation}
for all $\mathcal{H}^{N-1}\mres\dd^*E$-measurable sets $G\subset\Omega$. 
Let $\tau\in(0,1)$ and $B_r(x_0)\subset\Omega$. Without loss of generality, we may assume that $\tau<\frac 12$. We rescale $(E,u)$ in $B_1$ by setting $E_r=\frac{E-x_0}{r}$ and $u_r(y)=r^{-\frac 12}{u(x_0+ry)}$, for $y\in B_1$. We observe that $(E_r,u_r)$ satisfies the following $\Lambda r$-minimality relation: 
\begin{equation}
\label{eqq25}
\tilde{\mathcal{F}}(E_r,u_r):=\int_{B_1}\sigma_{E_r}|\D u_r|^2\,dy+\tilde{\Psi}_{E_r}(B_1)\leq \int_{B_1}\sigma_{E_r}|\D v|^2\,dy+\tilde{\Psi}_{F}(B_1)+\Lambda r|E_r\Delta F|,
\end{equation}
for any $(v,F)$ be such that $v-u\in H^1_0(B_1)$ and $E_r\Delta F\subset\subset B_1$. Here, we have denoted
\begin{equation*}
    \tilde{\Psi}_F(G)=\int_{\dd^*F\cap G}\theta(x_0+ry)\,d\mathcal{H}^{N-1}_y,\quad\forall G\subset\Omega.
\end{equation*}
We have to prove that there exists $\varepsilon_0=\varepsilon_0(\tau)$ such that, if $P(E_r;B_1)<\varepsilon_0$, then
\begin{equation}
\label{eqq26}
\tilde{\mathcal F}(E_r,u_r;B_{\tau })\leq C\big(\tau^N\tilde{\mathcal F}(E,u;B_1)+\tau^N r\big).
\end{equation}
For the rest of the proof, with a slight abuse of notation, we call $E_r$ by $E$ and $u_r$ by $u$. 
We note that, since $P(E;B_1)<\varepsilon_1$, by the relative isoperimetric inequality, either $|B_1\cap E|$ or $|B_1\setminus E|$ is small. Without loss of generality, We may assume that $|B_1\setminus E|\leq |B_1\cap E|$.
By the coarea formula, Chebyshev's inequality and the relative isoperimetric inequality, we may choose $\rho\in (\tau,2\tau)$ such that $\mathcal{H}^{n-1}(\partial^*E\cap \partial B_{\rho})=0$ and it holds
\begin{equation}
\label{eqq8}
\mathcal{H}^{n-1}(\partial B_\rho\setminus E)\leq \frac{C}{\tau} P(E;B_1)^{\frac{N}{N-1}}\leq \frac{C\varepsilon_0^{\frac{1}{N-1}}}{\tau} P(E;B_1),
\end{equation}
where $C=C(N)$. Now we test the minimality of $(u,E)$ with $(u,E\cup B_\rho)$.
We remark that, being $\mathcal{H}^{n-1}(\partial^*E\cap \partial B_{\rho})=0$, we can apply \cite[Proposition 2.2]{EspLamPis} with $U=F=B_1$ and $G=B_\rho$, thus obtaining
\begin{equation}
\label{eqq24}
 \tilde{\Psi}_{E\cup B_{\rho}}(B_1)
=\tilde{\Psi}_E(B_1\setminus \overline{B_\rho})+\tilde{\Psi}_{B_{\rho}}(B_1\setminus E^{(1)}).
 \end{equation}
Using the $\Lambda$-minimality relation \eqref{eqq25} with respect to the couple $(u,E\cup B_\rho)$, the equality \eqref{eqq24} to get rid of the common perimeter terms and recalling that $E=E^{(1)}$, we deduce
\begin{equation}
\int_{B_1}\sigma_E|\D u|^2\,dx+\tilde{\Psi}_E(B_{\rho})
\leq \int_{B_1}\sigma_{E\cup B_\rho}|\D u|^2\,dx+\tilde{\Psi}_{B_{\rho}}(B_1\setminus E^{(1)})+\Lambda r|B_\rho|.
\end{equation}
Taking into account the comparability to the perimeter \eqref{Comparability} and \eqref{eqq8}, recalling that $\rho\in (\tau,2\tau)$ and getting rid of the common Dirichlet terms, we deduce:
\begin{align}
\int_{B_{\tau}}\sigma_E|\D u|^2\,dx+K^{-1}P(E;B_{\tau})
& \leq \beta\int_{B_{2\tau}}|\D u|^2\,dx+K\mathcal{H}^{N-1}(\partial B_\rho\setminus E)+C(N,\Lambda) r\tau^N\\
&\leq \beta\int_{B_{2\tau}}|\D u|^2\,dx+\frac{C(N)K}{\tau} \varepsilon_1^{\frac{1}{N-1}}P(E;B_1)+C(N,\Lambda) r\tau^N.
\end{align}
Finally, we choose $\varepsilon_1$ such that

\begin{equation}
C(N)K \varepsilon_1^{\frac{1}{N-1}}\leq \tau^{N+1}\quad\mbox{and }\quad C(N)\varepsilon_1^{\frac{N}{N-1}} \leq \varepsilon_1(2\tau)|B_1|,
\end{equation}
where $\varepsilon_1$ corresponds to the $\varepsilon_0$ from \cite[Proposition 2.4]{FuscoJulin}, thus getting
\begin{equation}
\int_{B_{2\tau}}|\D u|^2\,dx\leq 2^n c_2\tau^n\int_{B_1}|\D u|^2\,dx.
\end{equation}
From this estimates \eqref{D} easily follows applying again the comparability to the perimeter \eqref{Comparability}.\\
\indent\textbf{Step 3:} \emph{Achieving the lower esitmate on the perimeter of $E$.} The proof matches exactly that of \cite[Proposition 4.4]{FuscoJulin}, given the comparability to the perimeter. We give only a sketch of the proof. We start by assuming that $x_0\in\displaystyle\dd^*E$. Without loss of generality, we may also assume that $x_0=0$. We denote by $C_1$ the constant $C$ appearing in \eqref{eqq11}, by $C_2$ the constant $C$ appearing in \eqref{D}. We recall that $\varepsilon_0$ is the constant appearing in Step 2. Arguing by contradiction, for $\tau\in\big(0,(2C_1)^{-2}\big)$ and $\sigma\in\big(0,\varepsilon_0(\tau)(2C_1 C_2)^{-1}\big)$ there exists a ball $B_r\subset U$, for some $r<\min\{\varepsilon_0(\tau),C_2\}$, such that 
\begin{equation*}
    P(E;B_r)\leq \varepsilon_0(\tau).
\end{equation*}
By using \eqref{eqq11} and \eqref{D}, we can easily prove by induction (see for example \cite[Theorem 4]{EspLam2} for the details) that
\begin{equation}
\label{Relazione iterativa}
\mathcal{F}(E,u;B_{\sigma\tau^hr})\leq\varepsilon_0(\tau)\tau^{\frac{h}{2}}(\sigma\tau^h r)^{N-1},\quad\forall h\in\N.
\end{equation}
From this estimate, we deduce that
\begin{equation*}
\lim_{\rho\rightarrow 0^+}\frac{P(E;B_\rho)}{\rho^{N-1}}=\lim_{h\rightarrow +\infty}\frac{P(E; B_{\sigma\tau^hr})}{(\sigma\tau^hr)^{N-1}}\leq\lim_{h\rightarrow+\infty}\varepsilon_0(\tau)\tau^{\frac{h}{2}}=0,
\end{equation*}
which implies that $x_0\not\in\dd^*E$, that is a contradiction.
If $x_0\in\dd E$, we get the same estimate by recalling that we chose the representative of $\dd E$ such that $\dd E=\overline{\dd\displaystyle^*E}$.\\
\indent\textbf{Step 4:} \emph{Proof of the Alhfors-regularity.} The proof of the final part of the statement follows as an application of the lower Ahlfors-regularity. Indeed, we have that
    \begin{equation*}
        \limsup_{r\rightarrow 0^+}\frac{\mathcal{H}^{N-1}(\dd^*E \cap B_r(x))}{r^{N-1}}=\limsup_{r\rightarrow 0^+}\frac{P(E;B_r(x))}{r^{N-1}}>0,\quad\forall x\in\dd E\cap \Omega.
    \end{equation*}
Thus, by \cite[(2.42)]{AFP}, we get $\mathcal{H}^{N-1}((\dd E\setminus \dd^*E)\cap\Omega)=0$. The Ahlfors-regularity of $\dd E$ follows as a consequence, taking also \eqref{DensityPerimeter} into account.

\end{proof}

The following results provide the main ingredients for the proof od Proposition \ref{rigot}, which in turn relies on the same strategy adopted in \cite[Lemma 3.6]{rigot}.

The subsequent lemma shows that optimal deformations satisfy a reverse H\"older inequality. Its proof is rather standard and it relies on the application of the Sobolev-Poincaré inequality.

\begin{Lem}\label{ReverseHolder}
    Let $(u,E)$ be a minimizer of \eqref{optimalDesign2}. There exists a positive constant $C_{0}=C_{0}\Big(N,\frac{\beta}{\alpha}\Big)$ such that it holds
\begin{equation}
\label{RevHol}
\int_{B_{\frac{r}{2}}(x_0)}|\D u|^2\,dx\leq C_{0}\bigg(\int_{B_r(x_0)}|\D u|^p\,dx\bigg)^{\frac{p}{2}},\quad\forall B_r(x_0)\subset\Omega,
\end{equation}
where $p=\frac{2N}{N+2}\in[1,2)$.
\end{Lem}
\begin{proof}
    Let $B_r(x_0)\subset\Omega$. Without loss of generality we may assume that $x_0=0$. Let $\xi\in C^1_c(B_r)$ be a cut-off function such that $\xi=1$ on $B_{\frac{r}{2}}$ and $|D\xi|\leq \frac{2}{r}$. By the Euler-Lagrange equation associated to $u$, taking as test function $(u-u_{B_r})\xi^2$ we get
\begin{align*}
    \int_{B_r} \sigma_E |\D u|^2\xi^2\,dx=-2\int_{B_r}\sigma_E \xi (u-u_{B_r})\prodscal{\D u}{\D \xi}\,dx.
\end{align*}
Applying H\"older's and Young's inequalities, we obtain
\begin{align*}
\int_{B_r}\sigma_E |\D u|^2\xi^2\,dx\leq \frac{1}{2}\int_{B_r}\sigma_E |\D u|^2\xi^2\,dx+2\int_{B_r}\sigma_E|u-u_{B_r}|^2|\D\xi|^2\,dx.
\end{align*}
Absorbing the first term in the right-hand side to the left-hand side and using that $\alpha\leq \sigma_E\leq \beta$, it holds that
\begin{align*}
\int_{B_{\frac{r}{2}}}|\D u|^2\,dx\leq \frac{1}{\alpha}\int_{B_r}\sigma_E |\D u|^2\xi^2\,dx\leq \frac{4}{\alpha} \int_{B_r}\sigma_E|u-u_{B_r}|^2|\D\xi|^2\,dx \leq \frac{16}{r^2}\frac{\beta}{\alpha}\int_{B_r}|u-u_{B_r}|^2\,dx.
\end{align*}
By the Sobolev-Poincaré inequality, we infer
\begin{align}
    \int_{B_{\frac{r}{2}}}|\D u|^2\,dx\leq C\bigg(\int_{B_r}|\D u|^p\,dx\bigg)^{\frac{p}{2}},
\end{align}
where $p=\frac{2N}{N+2}$ and $C=C\big(N,\frac{\beta}{\alpha}\big)$. Thus \eqref{RevHol} is proved.
\end{proof}

The following result establishes that it is possible to find a ball centered in a lower Ahlfors-regular curve where the rescaled Dirichlet energy is arbitrary small. Its proof relies on the strategy adopted in \cite[Corollary 38 in Section 23]{David}.

\begin{Lem}\label{LEM0}
Let $(u,E)$ be a minimizer of \eqref{optimalDesign2}, $U\subset\subset\Omega$ be an open set and $\Gamma\subset \partial E$ a subset which satisfies a  lower Ahlfors-regularity property with constant $C_\Gamma$. There exists $\varepsilon_0=\varepsilon_0\big(N,\frac{\alpha}{\beta}\big)>0$ such that for every $\varepsilon\in(0,\varepsilon_0)$ there exists a constant $a=a(N,K,\Lambda,C_\Gamma,\varepsilon)\in(0,1)$ such that for every $B_r(x_0)\subset U$, with $x_0\in\Gamma$ there exist $y\in\Gamma\cap B_{\frac{r}{3}}(x_0)$ and $t\in\big(ar,\frac{r}{3}\big)$ such that  
\begin{equation}
\int_{B_t(y)}|\D u|^2\,dx\leq \varepsilon t^{N-1}.
\end{equation}
\end{Lem}
\begin{proof} Let $\varepsilon>0$ and $B_r(x_0)\subset U$, with $x_0\in\dd E$. 
Let $p$ be the exponent in the reverse H\"older inequality given by Lemma \eqref{ReverseHolder} Setting
\begin{equation*}
\omega_p(x,s):=s^{1-\frac{2N}{p}}\bigg(\int_{B_s(x)}|\D u|^p\,dz\bigg)^{\frac{2}{p}},\quad\forall B_s(x)\subset\Omega,
\end{equation*}
we show that
\begin{equation}
\label{eqq6}
\int_{\Gamma\cap B_\frac{r}{3}(x_0)}\int_0^{\frac{r}{3}}\omega_p(y,t)\,\frac{d\mathcal{H}^{N-1}_y}{t}dt\leq C r^{N-1},
\end{equation}
for some positive contant $C=C(N,K,\Lambda,C_\Gamma)$.
Indeed, letting $\sigma\in[2,+\infty)$ be such that
\begin{equation}
\frac{1}{\sigma}+\frac{p}{2}=1,
\end{equation}
we choose $b\in\big(0,\frac{1}{\sigma}\big)$. By H\"older's inequality, we have
\begin{align}
\label{eqq16}
\omega_p(y,t)
& = t^{1-\frac{2N}{p}}\int_{B_t(y)}\bigg[|\D u(z)|^2\bigg(\frac{\mathrm{dist}(z,\Gamma)}{t}\bigg)^{\frac{2b}{p}}\bigg]^{\frac{p}{2}}\bigg(\frac{\mathrm{dist}(z,\Gamma)}{t}\bigg)^{-b}\,dz\\
& \leq t^{1-\frac{2N}{p}}\int_{B_t(y)}|\D u(z)|^2\bigg(\frac{\mathrm{dist}(z,\Gamma)}{t}\bigg)^{\frac{2b}{p}}\,dz\bigg[\int_{B_t(y)}\bigg(\frac{\mathrm{dist}(z,\Gamma)}{t}\bigg)^{-b\sigma}\,dz\bigg]^{\frac{2}{p\sigma}}.
\end{align}
We can partition $B_t(y)$ into subsets $B_k$ defined as
\begin{equation*}
B_k:=\{z\in B_t(y)\,:\, 2^{-k-1}t< \mathrm{dist}(z,\Gamma)\leq 2^{-k}t\},\quad\forall k\in\N_0,
\end{equation*}
and, since $\Gamma$ is Alhfors-regular, we apply \cite[Lemma 25]{David} to get
\begin{align}
\label{eqq15}
& \bigg[\int_{B_t(y)}\bigg(\frac{\mathrm{dist}(z,\Gamma)}{t}\bigg)^{-b\sigma}\,dz\bigg]^{\frac{2}{p\sigma}}\\
& =\sum_{k=0}^{+\infty} \int_{B_k}\bigg(\frac{\mathrm{dist}(z,\Gamma)}{t}\bigg)^{-b\sigma}\,dz\leq \sum_{k=0}^{+\infty} 2^{kb\sigma}|\{z\in B_t(y)\,:\,\mathrm{dist}(z,\Gamma)\leq 2^{-k}t\}|\leq Ct^N,
\end{align}
with $C=C(N,C_\Gamma)$. We remark that if $y\in\Gamma\cap B_{\frac{r}{3}}(x_0)$, $t\in\big(0,\frac{r}{3}\big)$, $z\in B_t(y)$, then
\begin{align}
z\in B_{\frac{2}{3}r}(x_0),\quad y\in\Gamma\cap B_t(z),\quad \mathrm{dist}(z,\Gamma)\in(0,t].
\end{align} 
Therefore, combining \eqref{eqq15} and \eqref{eqq16} and applying Fubini's theorem, the Ahlfors-regularity of $\Gamma$ and the upper bound on the energy (see \eqref{eqq11}), we get
\begin{align*}
& \int_{\Gamma \cap B_\frac{r}{3}(x_0)}\int_0^{\frac{r}{3}}\omega_p(y,t)\,\frac{d\mathcal{H}^{N-1}_y}{t}dt\\
& \leq C \int_{\Gamma \cap B_\frac{r}{3}(x_0)}\int_0^{\frac{r}{3}}t^{1-\frac{2N}{p}}t^{\frac{2N}{p\sigma}-1}\int_{B_t(y)}|\D u(z)|^2\bigg(\frac{\mathrm{dist}(z,\Gamma)}{t}\bigg)^{\frac{2b}{p}}\,dz\,d\mathcal{H}^{N-1}_ydt\\
& \leq C\int_{B_{\frac{2}{3}r}(x_0)}\int_{\mathrm{dist}(z,\Gamma)}^{\frac{r}{3}}t^{-N} \mathcal{H}^{N-1}(\Gamma\cap B_t(z))|\D u(z)|^2\bigg(\frac{\mathrm{dist}(z,\Gamma)}{t}\bigg)^{\frac{2b}{p}}\,dzdt\\
& \leq C\int_{B_{{\frac{2}{3}}r}(x_0)}|\D u(z)|^2\bigg[\int_{\mathrm{dist}(z,\Gamma)}^{\frac{r}{3}}\bigg(\frac{\mathrm{dist}(z,\Gamma)}{t}\bigg)^{\frac{2b}{p}}\bigg]\,dz\leq C r^{N-1},
\end{align*}
where $C=C(N,K,\Lambda,C_\Gamma)$, thus proving \eqref{eqq6}.\\
\indent We assume by contradiction that, for some $\varepsilon>0$ and for $a\in(0,1)$ to be chosen, for any $y\in \Gamma\cap B_{\frac{r}{3}}(x_0)$ and $t\in\big(a r,\frac{r}{3}\big)$, it holds that
\begin{equation}
\label{eqq7}
\omega_p(y,t)\geq \varepsilon^{\frac{2}{p}+1}.
\end{equation}
Thus, denoting
\begin{equation*}
\mathcal{A}:=(\Gamma\cap B_{\frac{r}{3}}(x_0))\times \bigg[ar,\frac{r}{3}\bigg],
\end{equation*}
by \eqref{eqq6}, \eqref{eqq7} and the Ahlfors-regularity of $\Gamma$, we obtain
\begin{align*}
r^{N-1} 
& \geq C\iint_{\mathcal{A}}\omega_p(y,t)\,\frac{d\mathcal{H}^{N-1}_y}{t}dt\geq C\varepsilon^{\frac{2}{p}+1}\mathcal{H}^{N-1}(\Gamma\cap B_{\frac{r}{3}}(x_0))\int_{ar}^{\frac{r}{3}}\,\frac{dt}{t}\\
& \geq C\varepsilon^{\frac{2}{p}+1} r^{N-1}\log\bigg(\frac{a}{3}\bigg),
\end{align*}
which is a contradiction if we choose $a=a(N,K,\Lambda,C_\Gamma,\varepsilon)>0$ sufficiently big. Therefore, there exist $y\in \Gamma\cap B_{\frac{r}{3}}(x_0)$ and $t\in\big(ar,\frac{r}{3}\big)$, where $a=a(N,K,\Lambda,C_\Gamma,\varepsilon)$ is a positive constant, such that $\omega_p(y,t)\leq \varepsilon^{\frac{2}{p}+1}$. By the reverse H\"older inequality \eqref{RevHol}, we conclude that
\begin{align}
\int_{B_{\frac{t}{2}}(y)}|\D u|^2\,dx
& \leq C\bigg(\int_{B_t(y)}|\D u|^p\,dx\bigg)^{\frac{p}{2}}= Ct^{N-\frac{p}{2}}\omega_p(y,t)^{\frac{p}{2}}\\
& \leq C \varepsilon\varepsilon^{\frac{p}{2}} t^{N-1}\leq \varepsilon t^{N-1},
\end{align}
where $C=C\big(N,\frac{\alpha}{\beta}\big)$,
provided that we choose $\varepsilon=\varepsilon\big(N,\frac{\alpha}{\beta}\big)>0$ sufficiently small.
\end{proof}

The following lemma establishes that if the rescaled Dirichlet integral is sufficiently small in some ball, then $E$ covers a large part of the ball. Combining this result with the previous lemma, we get a density estimate for $E$, which is given in Corollary \ref{COR1}.

\begin{Lem}
\label{LEM1}
Let $(u,E)$ be a minimizer of \eqref{optimalDesign2} and $U\subset\subset\Omega$ be an open set. There exists $\varepsilon_0=\varepsilon_0(\alpha,\beta,C_A)>0$ such that for every $\varepsilon\in(0,\varepsilon_0)$ 
there exist two positive constants $C_0=C_0(C_A)$ and $r_0=r_0(N,\Lambda,C_A)$ such that if for $B_r(x_0)\subset U$, with $r\in(0,r_0)$ and $x_0\in\dd E$, it holds
\begin{equation}
\label{eqq5}
\int_{B_r(x_0)}|\D u|^2\,dx\leq \varepsilon r^{N-1},
\end{equation}
then
\begin{equation}
\label{eqq4}
|B_r(x_0)\cap E|\geq C_0 r^N.
\end{equation}
\begin{proof}
Let $B_r(x_0)\subset U$ be such that
\begin{equation}
\label{eqq1}
\int_{B_r(x_0)}|\D u|^2\,dx\leq \varepsilon r^{N-1},
\end{equation}
where $\varepsilon>0$. We divide the proof in two steps.\\
\indent \textbf{Step 1:} We show that there exists a positive constant $r_0=r_0(N,\Lambda,C_A)$ such that
\begin{equation}
\label{eqq9}
\mathcal{H}^{N-1}(\dd B_r(x_0)\cap E)\geq \frac{C_A}{2}r^{N-1},
\end{equation}
for all $r\in(0,r_0)$.\\
\indent Applying the $\Lambda$-minimality of $(u,E)$ with respect to $(u,E\setminus B_r(x_0))$, we get
\begin{equation*}
\int_\Omega\sigma_E|\D u|^2\,dx+P(E;\Omega)\leq\int_{\Omega}\sigma_{E\setminus B_r(x_0)}|\D u|^2\,dx+P(E\setminus B_r(x_0);\Omega)+\Lambda|B_r(x_0)\cap E|.
\end{equation*}
Simplifying the previous inequality and estimating $|B_r(x_0)\cap E|\leq \omega_N r^N$, it follows that
\begin{equation*}
-(\beta-\alpha)\int_{B_r(x_0)\cap E}|\D u|^2\,dx+P(E;B_r(x_0))\leq \mathcal{H}^{N-1}(\dd B_r(x_0)\cap E)+\Lambda\omega_N r^N.
\end{equation*}
Thanks to \eqref{eqq1} and the Ahlfors regularity of $\dd E$, we get
\begin{equation*}
\big\{-(\beta-\alpha)\varepsilon+C_A-\Lambda\omega_N r\big\}r^{N-1}\leq \mathcal{H}^{N-1}(\dd B_r(x_0)\cap E).
\end{equation*}
Choosing $r_0=\frac{C_A}{4\Lambda\omega_N}$ and $\varepsilon_0<\frac{C_A}{4(\beta-\alpha)}$, we obtain \eqref{eqq9}.\\
\indent \textbf{Step 2:} We prove \eqref{eqq4}. Let us assume by contradiction that choosing $\delta<\frac{C_A}{16}$ and $M>0$ and $B_r(x_0)\subset U$, with $r<r_0$, \eqref{eqq5} holds and
\begin{equation}
\label{eqq3}
|B_r(x_0)\cap E|<\delta r^N,
\end{equation}
By Chebyshev inequality and Fubini's theorem, we get
\begin{align*}
&\mathcal{H}^1\bigg(\bigg\{\rho\in\bigg(\frac{r}{2},r\bigg)\,:\,\mathcal{H}^{N-1}(\dd B_\rho(x_0)\cap E)\geq\frac{4|B_r(x_0)\cap E|}{r}\bigg\}\bigg)\\
& \leq \frac{r}{4|B_r(x_0)\cap E|}\int_{\frac{r}{2}}^r\mathcal{H}^{N-1}(\dd B_\rho(x_0)\cap E)\,d\rho\leq \frac{r}{4}<\frac{r}{2}.
\end{align*}
This implies that there exists $\rho\in\big(\frac{r}{2},r\big)$ such that
\begin{equation}
\label{eqq2}
\mathcal{H}^{N-1}(\dd B_\rho(x_0)\cap E)\leq \frac{4|B_r(x_0)\cap E|}{r}\leq 4\delta r^{N-1}.
\end{equation}
Since $\delta<\frac{C_A}{16}$, we get a contradiction, \eqref{eqq9} being in force.
\end{proof}
\end{Lem}

\begin{Cor}
\label{COR1}
Let $(u,E)$ be a minimizer of \eqref{optimalDesign2} and $U\subset\subset\Omega$ be an open set. There exists a positive constant $C_0=C_0(N,\alpha,\beta,K,\Lambda,C_\Gamma)$ such that for every $B_{r}(x_0)\subset U$, with $x_0\in\dd E$, it holds 
\begin{equation}
|B_r(x_0)\cap E|\geq C_0 r^N.
\end{equation}
\end{Cor}
\begin{proof}
We fix $B_{r}(x_0)\subset U$, with $x_0\in\dd E$. Let us call by $\varepsilon_0$ the constant $\varepsilon_0$ appearing in Lemma \ref{LEM0}, by $\varepsilon_1$ and $r_0$ respectively the constants $\varepsilon_0$ and $r_0$ appearing in Lemma \ref{LEM1}. By Lemma \ref{LEM0} there exist a positive constant $a=a(N,K,\Lambda,C_\Gamma,\varepsilon)$ and $y\in\dd E\cap B_{\frac{r}{3}}(x_0)$ and $t\in\big(ar,\frac{r}{3}\big)$ such that  
\begin{equation}
\int_{B_t(y)}|\D u|^2\,dx\leq \varepsilon t^{N-1},
\end{equation}
where $\varepsilon\in(0,\min\{\varepsilon_0,\varepsilon_1\})$ and $r\in(0,r_0)$.
Applying Lemma \ref{LEM1}, we get that
\begin{equation*}
|B_r(x_0)\cap E|\geq |B_t(x_0)\cap E|\geq C t^N\geq C r^N,
\end{equation*}
for some positive constant $C=C(N,\alpha,\beta,K,\Lambda,C_\Gamma)$, which is the thesis.
\end{proof}

For what follows, it is useful to define the following function:
\begin{equation}
h(x,r)=r^{-N}\min\{|E\cap B_r(x)|,|B_r(x)\setminus E|\},\quad\forall B_r(x)\subset\Omega.
\end{equation}
In the following result, the boundary of an optimal shape $E$ is locally characterized in terms of the function $h$. It is the sets of points where the density of $E$ or $\Omega\setminus E$ is not too small.

\begin{Prop}
\label{PROP1}
Let $(u,E)$ be a minimizer of \eqref{optimalDesign2},  $U\subset\subset\Omega$ be an open set and $r_0>0$ be such that $\Lambda r_0\leq \frac{K^{-1}N}{2}$. There exists a positive constant $C_0=C_0(N,\alpha,\beta,K,\Lambda,C_\Gamma)$ such that it holds 
\begin{equation}
\dd E\cap U=\{x\in U\,:\, h(x,r)\geq C_0,\,\forall r\in(0,\min\{r_0,\mathrm{dist}(x,\dd U)\})\}.
\end{equation}
\end{Prop}
\begin{proof}
Let $x_0\in\dd E\cap U$ and $r\in(0,\min\big\{r_0,\mathrm{dist}(x_0,\dd U),\frac{K^{-1}N}{2}\big\})$. If $h(x_0,r)=r^{-N}|B_r(x_0)\cap E|$, by Corollary \ref{COR1}, there exists a positive constant $C_1=C_1(N,\alpha,\beta,K,\Lambda,C_\Gamma)$ such that
\begin{equation}
|B_r(x_0)\cap E|\geq C_1 r^N,
\end{equation}
that is $h(x_0,r)\geq C_1$. If $h(x_0,r)=r^{-N}|B_r(x_0)\setminus E|$, we assume that $\mathcal{H}^{N-1}(\dd^*E\cap\dd B_r(x_0))=0$. We fix $s\in(r,\min\{r_0,\mathrm{dist}(x_0,\dd U)\})$.
Let $F\subset\R^N$ be such that $E\Delta F\subset\subset B_s(x_0)$. By the isoperimetric inequality we infer that
\begin{align*}
    |E\Delta F|
    & =|E\Delta F|^{\frac{N-1}{N}}|E\Delta F|^{\frac{1}{N}}\leq \frac{P(E\Delta F)}{N}s=\frac{s}{N}P(E\Delta F;B_s(x_0))\\
    & \leq \frac{s}{N}[P(E;B_s(x_0))+P(F;B_s(x_0))].
\end{align*}
Thus, taking \eqref{Comparability} into account, we get
\begin{equation}
    (\beta-\alpha)\int_{B_s(x_0)\setminus E}|\D u|^2\,dx+\bigg(K^{-1}-\frac{\Lambda s}{N}\bigg)P(E;B_s(x_0))\leq\bigg(K+\frac{\Lambda s}{N}\bigg)P(F;B_s(x_0)).
\end{equation}
Testing the previous relation with the set $E\cup B_r(x_0)$ and letting $s\rightarrow r^+$, we get
\begin{equation}
(\beta-\alpha)\int_{B_r(x_0)\setminus E}|\D u|^2\,dx+\bigg(K^{-1}-\frac{\Lambda r}{N}\bigg)P(E;B_r(x_0))\leq\bigg(K+\frac{\Lambda r}{N}\bigg)\mathcal{H}^{N-1}(\dd B_r(x_0)\setminus E).
\end{equation}
Thus, adding $\big(K^{-1}-\frac{\Lambda r}{N}\big)\mathcal{H}^{N-1}(\dd B_r(x_0)\setminus E)$ to both sides of the previous inequality, we get
\begin{equation}
\bigg(K^{-1}-\frac{\Lambda r}{N}\bigg) P(B_r(x_0)\setminus E)\leq \big(K+K^{-1}\big)\mathcal{H}^{N-1}(\dd B_r(x_0)\setminus E).
\end{equation}
Taking into account that $1-\frac{\Lambda r_1}{N}\geq \frac{K^{-1}N}{2}$, the isoperimetric inequality in the left hand-side yields
\begin{equation}
|B_r(x_0)\setminus E|^{\frac{N-1}{N}}\leq c(N,K)\mathcal{H}^{N-1}(\dd B_r(x_0)\setminus E).
\end{equation}
Setting $m(r):=|B_r(x_0)\setminus E|$, the previous inequality can be rephrased as
\begin{equation}
m(r)^{\frac{N-1}{N}}\leq c(N,K) m'(r),\quad\text{for a.e. }r\in(0,\min\{r_0,\mathrm{dist}(x_0,\dd U)\}).
\end{equation}
Integrating the previous inequality in $(0,r)$, we get
\begin{equation}
m(r)\geq c(N,K)r^N,\quad\forall r\in(0,\min\{r_0,\mathrm{dist}(x_0,\dd U)\}),
\end{equation}
which implies that $h(x_0,r)\geq c(N,K)$.
\end{proof}


The main result of this section follows.

\begin{Prop}\label{rigot}
Let $(u,E)$ be a minimizer of \eqref{optimalDesign2} and $U\subset\subset\Omega$ be an open set. Then the open sets
\begin{equation}
E_0=\{x\in U\,:\,\exists r>0\text{ s.t. }|B_r(x)\setminus E|=0\},
\end{equation}
\begin{equation}
E_1=\{x\in U\,:\,\exists r>0\text{ s.t. }|B_r(x)\cap E|=0\}
\end{equation}
satisfy the condition-$B$ with some positive constants $C_0$ and $r_0$, and $E_0$ is equivalent to $E$. Moreover, it holds that $E_1=\R^N\setminus\overline{E_0}$ and $\dd E_0=\dd E_1=\dd E$.  
\end{Prop}
\begin{proof}
We only prove that $E_0$ and $E_1$ satisfy the condition-$B$. Indeed, the validity of the other statements has been showed in \cite[Lemma 3.6]{rigot}. Setting $r_0:=\frac{1}{2}\min\big\{\frac{K^{-1}N}{2\Lambda},\mathrm{dist}(U,\dd\Omega)\big\}$, let $x_0\in\dd E$ and $r\in(0,r_0)$. 
By Proposition \ref{PROP1}, there exists a positive constant $C_1=C_1(N,\alpha,\beta,K,\Lambda,C_A)$ such that $h\big(x_0,\frac{r}{2}\big)\geq C_1$. Thus
\begin{equation}
\label{eqq10}
|E_0\cap B_{\frac{r}{2}}(x_0)|\geq 2^{-N}C_1 r^N\quad\text{and}\quad |E_1\cap B_{\frac{r}{2}}(x_0)|\geq 2^{-N}C_1 r^N.
\end{equation}
Since, by the Ahlfors regularity of $\dd E$, it holds that 
\begin{equation*}
\bigg|\bigg\{y\in B_{\frac{r}{2}}(x_0)\,:\,\mathrm{dist}(y,\dd E)\leq\frac{tr}{2}\bigg\}\bigg|\leq Ctr^N,\quad\forall,t\in(0,1),
\end{equation*}
where $C=C(N,C_A)>0$
(see \cite[Lemma 3.6]{rigot} or \cite[Lemma 25 in Chapter 23]{David}). Choosing $t\in\big(0,\frac{C_3}{2^{N+1}C(N,C_A)}\big)$, by \eqref{eqq10} we get that
\begin{equation}
\bigg|\bigg\{y\in B_{\frac{r}{2}}(x_0)\,:\,\mathrm{dist}(y,\dd E)\leq\frac{tr}{2}\bigg\}\bigg|< \min\{|E_0\cap B_{\frac{r}{2}}(x_0)|,|E_1\cap B_{\frac{r}{2}}(x_0)|\}.
\end{equation}
As a consequence, there exist $y_0\in E_0\cap B_{\frac{r}{2}}(x_0)$ and $y_1\in E_1\cap B_{\frac{r}{2}}(x_0)$ such that 
\begin{equation}
\mathrm{dist}(y_0,\dd E)>\frac{tr}{2}\quad\text{and}\quad \mathrm{dist}(y_1,\dd E)>\frac{tr}{2},
\end{equation}
which is our aim.
\end{proof}

As an application to condition-B, we can give an alternative proof for an estimate that was already found by Larsen, stated in \cite[Theorem 4.1]{Larsen2003} and proved with completely different arguments.  

\begin{Prop}\label{Larsen4}
    Let $(u,E)$ be a minimizer of $\mathcal{F}$ defined in \eqref{optimalDesign2} and let $A_i$ be a connected component of $E$. Then there exists two positive constants $C_1$ $C_2$ such that if $|A_i|\leq C_1$, then
    \begin{equation*}
        \Big|\Big\{x\in E\setminus A_i\,:\,\mathrm{dist}(x,A_i)\leq C_2|A_i|^{\frac{1}{N}}\Big\}\Big|>|A_i|.
    \end{equation*}
\end{Prop}

    \begin{proof}
        Let $x\in A_i$. Since $E$ satisfies condition-B with some constants $C_0$ and $r_0$, setting $r:=\Big(\frac{2|A_i|}{\omega_N}\Big)^{\frac{1}{N}}\leq C_0 r_0$, there exists $y\in E$ such that
        \begin{equation*}
            B_r(y)\subset B(x,C_0r)\cap E.
        \end{equation*}
        By the choice of $r$, we get that $|B_r(y)|>|A_i|$. Thus, it holds necessarily that
        \begin{equation*}
    B_r(y)\subset E\setminus A_i.
        \end{equation*}
        Furthermore, for any $z\in B_r(y)$, we have that
    \begin{equation*}
        \mathrm{dist}(z,A_i)\leq \frac{r}{C_0}=\frac{1}{C_0}\bigg(\frac{2}{\omega_N}\bigg)^\frac{1}{N}|A_i|^{\frac{1}{N}}.
    \end{equation*}
    Accordingly, 
    \begin{align}
        \Bigg|\Bigg\{z\in E\setminus A_i\,:\,\mathrm{dist}(z,A_i)\leq \frac{1}{C_0}\bigg(\frac{2}{\omega_N}\bigg)^\frac{1}{N}|A_i|^{\frac{1}{N}}\Bigg\}\Bigg|\geq |B_r(y)|>|A_i|,
    \end{align}
    which is the thesis.
    \end{proof}


\section{Control of the flatness in a quantitative way}


Lemma \ref{variantB} below is a variant of the well-known classical fact that for a uniformly rectifiable set, one can find many balls in which the bilateral flatness is small. In other words, the set of points where the flatness is larger than any threshold $\varepsilon>0$, is a porous set. The following variant says that one can even choose the center of the ball, where the set is flat, in a given subset $\Gamma \subset \partial E$ provided that $\Gamma$ satisfies itself a lower Ahlfors regularity estimate. 

\begin{Def}\label{lowerA} We say that $\Gamma\subset \R^N$ satisfies a lower Ahlfors-regularity estimate if there exists $r_0,C_\Gamma>0$ such that for all $r\in(0,r_0)$ and  $x\in \Gamma$, it holds
$$\mathcal{H}^{N-1}(\Gamma\cap B_r(x))\geq C_\Gamma r^{N-1}.$$  
\end{Def}

Here below is the lemma that will be needed to control the flatness in many balls.

\begin{Lem}\label{variantB} Let $(u,E)$ be a minimizer of $\mathcal{F}$ defined in \eqref{optimalDesign22} and let $\Gamma \subset \partial E$ be a subset that satisfies a lower Ahlfors-regularity estimate with constants $r_0, C_\Gamma$ (see Definition \ref{lowerA}). Then,   for every $\varepsilon>0$ there exists $a>0$ such that the following holds: for all $x \in \Gamma$ and $r\in (0,r_0/2)$, there exist $y\in \Gamma \cap B_r(x)$ and $t\in (ar,r)$ such that 
$$\beta(y,t)\leq \varepsilon.$$
\end{Lem}
\begin{proof} Since $\partial E$ is Ahlfors-regular and uniformly rectifiable, we know by Proposition \ref{EquivBWGL} that $\partial E$ satisfies the bilateral weak geometric lemma, i.e. for every $\varepsilon>0$, the set
\begin{equation}
A:=\{(x,t) \in \partial E\times \R_+ \; :\; \beta(x,t)>\varepsilon\} \label{definitionA}
\end{equation}
is a Carleson set (see \cite[Definition 2.2]{david-semmes}). More precisely there exists $C>0$ and  such that 
\begin{equation}
\int_0^r \int_{\partial E\cap B_r(x)} {\bf 1}_{A}(y,t) d\mathcal{H}^{N-1}(y)\frac{dt}{t}\leq Cr^{N-1}, \label{acontredire}
\end{equation}
for all $r>0$. 

Now, let $\Gamma\subset \partial E$ satisfy a lower Ahlfors regularity estimate with constants $r_0,C_\Gamma>0$ and let $\varepsilon>0$ be given. Fix also some $x\in \Gamma$ and $0<r\leq r_0/2$. Assume by contradiction that for all  $t\in (ar,r)$ and for all $y \in \Gamma \cap B_r(x)$ it holds
\begin{equation}
\beta(y,t)> \varepsilon. \label{autiliser}
\end{equation}
Defining  
$$A_{\Gamma}:=\{(x,t) \in \Gamma\times \R_+ \; :\; \beta(x,t)>\varepsilon\},$$
using moreover that  \eqref{autiliser} holds for all $t\in (ar,r)$, we get
\begin{eqnarray}
\int_0^r \int_{\partial E\cap B_r(x)} {\bf 1}_{A}(y,t) \,d\mathcal{H}^{N-1}(y)\frac{dt}{t}&\geq& \int_0^r \int_{\partial E\cap B_r(x)} {\bf 1}_{A_{\Gamma}}(y,t) \,d\mathcal{H}^{N-1}(y)\frac{dt}{t}, \notag \\
&\geq & \int_{ar}^r \int_{\partial E\cap B_r(x)} {\bf 1}_{A_{\Gamma}}(y,t) \,d\mathcal{H}^{N-1}(y)\frac{dt}{t}, \notag \\
&\geq & \int_{ar}^r  \mathcal{H}^{N-1}(\Gamma \cap B_r(x))\,\frac{dt}{t}, \notag \\
&\geq & C_\Gamma r^{N-1} \ln(1/a), \notag 
\end{eqnarray}
which contradicts \eqref{acontredire} for $a$ small enough, and the lemma follows. 
\end{proof}

\begin{Lem}
\label{LemmaFlatnessExcess}
For all $\varepsilon_0>0$ there exists $\eta_0>0$ such that the following holds: for any minimizer of the functional $\mathcal{F}$ defined in \eqref{optimalDesign22} in $B_r(x)$ satisfying 
$$\beta(x,r)\leq \eta_0,$$
we have
$$\mathbf{e}(x,r/2)\leq \varepsilon_0.$$
\end{Lem}
\begin{proof}
Let us assume by contradiction that there exists $\varepsilon_0>0$, a sequence of balls $B_{r_h}(x_h)$ and a sequence of $\Lambda$-minimizers $(E_h,u_h)$ of $\mathcal{F}$ in $B_{r_h}(x_h)$ such that 
\begin{equation*}
    \beta_h(x_h,r_h)  \rightarrow 0,
\end{equation*}
as $h\rightarrow +\infty$, and 
\begin{equation*}
    \mathbf{e}_h(x_h,r_h)> \varepsilon_0,\quad\forall h\in\N.
\end{equation*}
Setting $(E'_h,u'_h)=\Big(r_h^{-1}(E_h-x_h), r_h^{-\frac{1}{2}}u_h(x_h+r_h\cdot)\Big)$, it holds that $(E'_h,u_h)$ is a $\Lambda r_h$-minimizer of $\mathcal{F}$ in $B_1$. Furthermore, by scaling, it follows that
\begin{equation}
    \beta_h(0,1)  \rightarrow 0, \label{convergenceBeta}
\end{equation}
as $h\rightarrow +\infty$, and 
\begin{equation}
\label{eqq13}
    \mathbf{e}_h(0,1)> \varepsilon_0,\quad\forall h\in\N.
\end{equation}
Since $E_h$ is uniformly Ahlfors-regular, we know that 
$$P(E_h;B_1)\leq C_A, \quad \forall h\in\N,$$
where $C_A$ is a constant independent of $h$. This allows to extract a further subsequence such that $E_h\to E$ in $L^1(B_1)$ and,  from \eqref{convergenceBeta}, we infer that, up to a rotation, $E=B_1\cap \{x_N>0\}$. Let $C:={\mathbf C}(0,1/2,e_N)$ be a cylinder satisfying $B_{\frac{1}{2}}\subset C \subset \mathbf{C}_1$. Then by \cite[Proposition 22.6]{maggi} we know that 
$$\mathbf{e}_h(0,1/2)\to 0,$$
which contradicts \eqref{eqq13}, and concludes the proof.
\end{proof}



\begin{Lem}\label{porosity} There exists $\varepsilon_0>0$ such that for any $\varepsilon\in(0,\varepsilon_0)$ there exists $a>0$ and $r_0>0$ such that the following holds. Let  $(u,E)$ be a minimizer for of $\mathcal{F}$ defined in \eqref{optimalDesign22} and    $\Gamma\subset \partial E$ a subset which satisfies a lower Ahlfors-regularity property with constant $C_\Gamma$.  Then, for any $x\in \Gamma$ and $r\in(0,r_0)$ there exists $y\in \Gamma \cap B_r(x)$ such that 
$$e(y,ar) +\omega(y,ar) + \Lambda ar  \leq \varepsilon_0.$$
\end{Lem}

\begin{proof} 
Let us denote by $\varepsilon_1$ the constant $\varepsilon_0$ that appears in Lemma \ref{LEM0} and let $\eta_0$ be the constant of Lemma \ref{LemmaFlatnessExcess}. We choose $\varepsilon<3\varepsilon_1$. By Lemma \ref{LEM1}, there exists $b>0$ such that for any $x_0\in\Gamma$ and $B_\rho(x_0)\subset U$ there exist $z\in B_{\frac{\rho}{3}}(x_0)$ and $s\in\big(b\rho,\frac{\rho}{3}\big)$ such that
\begin{equation}
\label{eqq18}
\omega(z,s)\leq\frac{\varepsilon}{3}.
\end{equation}
Furthermore, by Lemma \ref{LemmaFlatnessExcess}, there exists $\eta_0>0$ such that if, for any $B_r(x_0)\subset\Omega$ such that $\beta(x_0,r)<\eta_0$, we have
\begin{align}
\label{eqq17}
\mathbf{e}\bigg(x_0,\frac{t}{2}\bigg)\leq b^{N-1}\frac{\varepsilon}{3}
\end{align}
Finally, by Lemma \ref{variantB}, there exists $a>0$ such that for any $x_0\in\Gamma$ and $r\in\big(0,\frac{r_0}{2}\big)$ there exist $y\in\Gamma\cap B_r(x)$ and $t\in(a r,r)$ such that
\begin{equation*}
    \beta(y,t)\leq \eta_0.
\end{equation*}
Thus, we get that 
\begin{equation*}
    \mathbf{e}\bigg(y,\frac{r}{2}\bigg)\leq b^{N-1}\frac{\varepsilon}{3}.
\end{equation*}
 We choose $\rho=\frac{t}{2}$. By the previous inequality, it follows that
 \begin{equation*}
     \mathbf{e}(z,s)\leq \bigg(\frac{t}{2s}\bigg)^{N-1}\mathbf{e}\bigg(y,\frac{t}{2}\bigg)\leq \frac{\varepsilon}{3}.
 \end{equation*}
If we chose $r_0<\frac{2\varepsilon}{\Lambda}$, then, by the previous chain of inequalities and \eqref{eqq18}, we get
\begin{equation*}
    \mathbf{e}(z,s)+\omega(z,s)+\Lambda s<\varepsilon,
\end{equation*}
which is the thesis.
\end{proof}

\begin{Cor}\label{hausdorff} $\partial E_{sing}$ has Hausdorff dimension stricitly less than 
$N-1$.
\end{Cor} 
\begin{proof}  By the previous lemma, $\partial E_{sing}$  is a porous set in $\partial E$. Then it is standard to show that the Hausdorff dimension is strictly less than $N-1$ (see for instance \cite[Remark 3.29]{rigot}).
\end{proof}


\section{Regularity in dimension $N=2$}

Throughout the section, we derive some properties of solutions of the optimal design problem \eqref{optimalDesign} in the case $N=2$. In Subsection \ref{Monoton}, we employ an argument due to Bonnet \cite{bonnet} to show the validity of a monotonicity formula concerning $\D u$ (see Proposition \ref{monotonie}), provided that the ratio between the coefficients $\alpha$ and $\beta$ is smaller than 4. The smoothness of the free boundary, which is the statement in Theorem \ref{main1}, follows as a consequence. In Subsection \ref{TwoPoints}, we study the spectral problem associated to the optimal design when the set $E$ meets a circle exactly in two points (see Proposition \eqref{spectral1}). We use the estimate on the first eigenvalue to derive a partial regularity result concerning the free boundary, which is contained in Corollary \ref{regularity2phases}. Finally, in Subsection \ref{SectionDist}, we prove a quantitative estimate of the mutual distance between connected components of the optimal set.


    \subsection{Full regularity in the case $\beta < 4\alpha$ from a  monotonicity formula}
    \label{Monoton}

  \label{sectionsmooth}

We will need in the sequel the following lemma which is the justification of an integration by parts formula.

 \begin{Lem}[Integration by parts]\label{IPP}
 Let    $\alpha \leq \beta$ and  $\sigma\colon\Omega \to \R^+$ any given measurable function such that $\alpha \leq \sigma \leq \beta$. Let   $u$ be a local  minimizer for the energy
 $$\int_{\Omega} \sigma |\nabla u|^2 \;dx.$$
 Then for all $x_0\in \Omega$ and a.e. $r\in(0, \mathrm{dist}(x_0,\partial \Omega))$ we have
 $$\int_{B_r(x_0)}\sigma |\nabla u|^2 \,dx= \int_{\partial B_r(x_0)}\sigma u \frac{\partial u}{\partial \nu}\,dx,$$
 and
 $$  \int_{\partial B_r(x_0)}\sigma   \frac{\partial u}{\partial \nu}\;dx = 0.$$
  \end{Lem}  
  \begin{proof} Without loss of generality, we may assume that $x_0=0$. We know from the local minimality of $u$ that it holds
  \begin{equation}
  \int_{\Omega} \sigma \nabla u \cdot \nabla \varphi \,dx = 0,  \quad\forall  \varphi \in H^1_0(\Omega).
  \label{opti}
  \end{equation}
  For $\varepsilon>0$, let us choose $\varphi= \varphi_\varepsilon u$ with $\varphi_\varepsilon:=g_\varepsilon(|x|)$ where $g_\varepsilon\colon[0,+\infty)\rightarrow\R$ is the continuous function defined as follows: $g_\varepsilon(t)=1$ for any $t\in[0,(1-\varepsilon)r]$, $g_\varepsilon(t)=0$ for any $t\in[(1+\varepsilon)r,+\infty[$ and it is linear on $[(1-\varepsilon)r,(1+\varepsilon)r]$. Applying \eqref{opti}, we get 
  $$\int_{\Omega} \sigma \nabla u \cdot u\nabla \varphi_\varepsilon \,dx + \int_{\Omega} \sigma \nabla u \cdot \varphi_\varepsilon \nabla u \,dx= 0 .$$
  It is clear that $\varphi_\varepsilon$ converges strongly in $L^2$ to $\mathbbm{1}_{B_r}$, which implies
  $$ \lim_{\varepsilon\rightarrow 0}\int_{\Omega} \sigma \nabla u \cdot \varphi_\varepsilon r \nabla u \,dx =\int_{B_r} \sigma |\nabla u |^2 \,dx.$$
  On the other hand $\varphi_\varepsilon$ is Lipschitz continuous and $\D\varphi_\varepsilon(x)=\frac{x}{2\varepsilon r|x|}\mathbbm{1}_{B_{(1+\varepsilon)r}\setminus B_{(1-\varepsilon)r}}(x)$ for almost every $x\in\Omega$, so that 
  $$\int_{\Omega} \sigma \nabla u \cdot u\nabla \varphi_\varepsilon \,dx=\frac{1}{2\varepsilon r}\int_{B_{(1+\varepsilon)r}\setminus B_{(1-\varepsilon)r}} \sigma u \nabla u \cdot \frac{x}{|x|} \,dx$$
which converges to $ \int_{\partial B_r} \sigma u \nabla u \cdot \nu \;d\mathcal{H}^1$ for a.e. $r$ by the Lebesgue's differentiation theorem. Therefore, passing to the limit we obtain 
 $$\int_{B_r}\sigma |\nabla u|^2 \,dx= \int_{\partial B_r}\sigma u \frac{\partial u}{\partial \nu}\,dx,$$
 which proves the first half of the statement. For the second assertion, we take $\phi=\varphi_\varepsilon$ itself as a test function. This gives 
 $$\int_{\Omega} \sigma \nabla u \cdot \nabla \varphi_\varepsilon \;dx  = 0 ,$$
 which actually yields
 $$\frac{1}{2\varepsilon r}\int_{{B_{(1+\varepsilon)r}\setminus B_{(1-\varepsilon)r}}} \sigma  \nabla u \cdot \frac{x}{|x|} \,dx=0.$$
 The lemma follows from passing to the limit as $\varepsilon \to 0$.
\end{proof}

We are now able to prove the monotonicity formula.
  
    \begin{Prop}\label{monotonie}
 Let  $\sigma\colon\Omega \to \R^+$ a given function such that 
 $$\alpha \leq \sigma \leq \beta.$$
Let   $u$ be a local  minimizer for the energy
 $$\int_{\Omega} \sigma |\nabla u|^2 \;dx.$$ 
 Then, for all $x_0\in \Omega$ and $r\in(0, \mathrm{dist}(x_0,\partial \Omega))$ the function  
 $$r\mapsto  \frac{1}{r^{\gamma}} \int_{B_r(x_0)} \sigma |\nabla u|^2 \,dx$$
 is nondecreasing with $\gamma=2\sqrt{\frac{\alpha}{\beta}}$.
  \end{Prop}  

\begin{proof}  Without loss of generality, we may assume that $x_0=0$. We define 
$$E(r):=\int_{B_r}\sigma |\nabla u|^2 \;dx.$$
Notice that  $E(r)$ can be rewritten in the form
$$E(r)=\int_{0}^r \int_{\partial B_t} \sigma| \nabla u|^2 \, d\mathcal{H}^1 dt,$$
which directly shows that $E(r)$ is absolutely continuous and, in particular, differentiable   almost everywhere, with
$$E'(r)=\int_{\partial B_r} \sigma| \nabla u|^2 \, d\mathcal{H}^1.$$
We recall the following classical Wirtinger inequality, valid for any function $u \in W^{1,2}(\partial B_r)$,
\begin{eqnarray}
\int_{\partial B_r}  (u-m_r)^2 \,dx \leq r^2  \int_{\partial B_r} \left( \frac{\partial u}{ \partial \tau}\right)^2 \,dx,
\end{eqnarray}
where $m_r$ is the average of $u$, namely, $m_r:=\frac{1}{2\pi r}\int_{\partial B_r} u \; d\mathcal{H}^1$.
From the minimizing property of $u$ we may apply Lemma \ref{IPP} which yields
$$\int_{B_r}\sigma |\nabla u|^2 \,dx= \int_{\partial B_r}\sigma u \frac{\partial u}{\partial \nu}\,dx,\quad \text{for a.e. }r>0,$$
and
$$  \int_{\partial B_r}\sigma   \frac{\partial u}{\partial \nu}\;dx = 0.$$
Using the equations above, H\"older's and Young's  inequalities, we can write
\begin{eqnarray}
E(r)=\int_{B_r}\sigma |\nabla u|^2 \,dx &=& \int_{\partial B_r}\sigma u \frac{\partial u}{\partial \nu}\,dx = \int_{\partial B_r}\sigma(u-m_r) \frac{\partial u}{\partial \nu}\,dx \notag \\
&\leq &   \left( \int_{\partial B_r}\sigma(u-m_r)^2 \bigg( \frac{\partial u}{\partial \nu}\bigg)^2\,dx \right)^{1/2} \notag \\
&\leq & \frac{\varepsilon}{2} \int_{\partial B_r}\sigma(u-m_r)^2\,dx  + \frac{1}{2\varepsilon}\int_{\partial B_r}\sigma\bigg( \frac{\partial u}{\partial \nu}\bigg)^2\,dx   \notag \\
&\leq &  \frac{\varepsilon}{2}\beta \int_{\partial B_r} (u-m_r)^2  \,dx+ \frac{1}{2\varepsilon}\int_{\partial B_r}\sigma\bigg( \frac{\partial u}{\partial \nu}\bigg)^2\,dx \notag \\
&\leq & \frac{\varepsilon}{2}\beta  r^2  \int_{\partial B_r} \left( \frac{\partial u}{ \partial \tau}\right)^2 \,dx+ \frac{1}{2\varepsilon}\int_{\partial B_r}\sigma\bigg( \frac{\partial u}{\partial \nu}\bigg)^2\,dx \notag  \\
&\leq & \frac{\varepsilon}{2} \frac{\beta}{\alpha}  r^2  \int_{\partial B_r} \sigma \left( \frac{\partial u}{ \partial \tau}\right)^2 \,dx+ \frac{1}{2\varepsilon}\int_{\partial B_r}\sigma\bigg( \frac{\partial u}{\partial \nu}\bigg)^2\,dx,\notag 
\end{eqnarray}
where $\varepsilon>0$.
Finally, by choosing $\varepsilon=\frac{1}{r}\sqrt{\frac{\alpha}{\beta}}$ we find that
$$E(r)\leq \frac{r}{2}\sqrt{\frac{\beta}{\alpha}}E'(r).$$
This implies that the function
 $$r\mapsto  \frac{1}{r^{\gamma}} \int_{B_r(x)} \sigma |\nabla u|^2 \;dx$$
 is nondecreasing with $\gamma=2\sqrt{\frac{\alpha}{\beta}}$.
\end{proof}

As a consequence, we obtain the result already announced in the introduction as  Theorem~\ref{main1}. 


\begin{proof}[Proof of Theorem \ref{main1}] Let $(u,E)$ be a minimizer for Problem \eqref{optimalDesign} with  $\beta< 4\alpha$. Then by \cite[Theorem 1]{espFusco}, we know that $(u,E)$ is a $\Lambda$-minimizer, i.e. $(u,E)$ is a minimizer of the functional
\begin{equation}
 \int_{\Omega} \sigma_{E} |\nabla u|^2 \;dx + P(E; \Omega) + \Lambda||E|-V_0|. 
\end{equation}

In particular, the set $E$ being fixed, we infer that $u$ must be  a minimizer of the energy $\int_{\Omega} \sigma_{E} |\nabla u|^2 \;dx$,
and according to Proposition \ref{monotonie},
we know that there exists $r_0>0$ such that for all $x_0\in \partial E$ and $r\in (0,r_0)$ we have 
$$\int_{B_r(x_0)} \sigma_E |\nabla u|^2 \;dx \leq C_0r^{1+\varepsilon},$$
with $\varepsilon=2\sqrt{\alpha/\beta}-1>0$ and $C_0=r_0^{-2\sqrt{\alpha/\beta}}\int_{B_{r_0}(x_0)} \sigma_E |\nabla u|^2 \;dx $.

Now let  $F\subset\R^N$ be any set of finite perimeter such that $E\Delta F\subset B_r(x_0)$. Then, testing the minimality of $(u,E)$ with the competitor $(u,F)$ we get
\begin{equation}
 \int_{\Omega} \sigma_{E} |\nabla u|^2 \,dx + P(E; \Omega)\leq \int_{\Omega} \sigma_{F} |\nabla u|^2 \,dx + P(F; \Omega) + \Lambda||F|-V_0|,
\end{equation}
which implies, since $E\Delta F\subset B_r(x_0)$,
\begin{align}
  P(E; B_r(x_0))   &\leq P(F; B_r(x_0)) + C \int_{B_r(x_0)} \sigma_{E} |\nabla u|^2 \,dx  + C|B_r(x_0)|, \notag \\
  &\leq  P(F; B_r(x_0)) + C r^{1+\varepsilon}.
\end{align}

 This means that $E$ falls into the theory of almost minimizers for the perimeter, and then by the classical result of Tamanini \cite{tamanini} we deduce that the singular set is regular up to a singular set of dimension $N-8$. Since here $N=2$, the singular set is actually empty.
\end{proof}

 \subsection{Monotonicity formula for a boundary intersecting by only two points}
 \label{TwoPoints}
 
In  the previous subsection, we have obtained a monotonicity behavior of the energy by using the classical Wirtinger inequality in order to estimate the derivative of the energy with respect to the radius of the ball. This strategy is necessarily non optimal, due to the coefficients $\sigma$ in front of the energy.

In this section, we try to improve the monotonicity behavior of the energy by analysing  precisely the  Wirtinger constant  taking into account the weight $\sigma$ in the inequality. In other words, we arrive to a new spectral problem on the circle, with weight $\sigma$.

In the case when the two regions $\{\sigma=\alpha\}$ and $\{\sigma=\beta\}$ are both connected on the circle, we obtain a good decay behavior of the energy of the type $o(r)$ leading to $C^1$-regularity (see Corollary \ref{regularity2phases}).   It was surprising to the authors that even in such ``easy'' case, the associated 1D-spectral problem on the circle was so difficult to compute (see Proposition \ref{spectral1}).

Moreoever when  those regions are not connected, the computations become really painful and some numerical evidences shows that one cannot hope to obtain a good behavior of the energy in full generality. In other words, using this strategy it seems difficult to prove a full regularity result without any restriction on $\alpha$, $\beta$ or the regions $\{\sigma=\alpha\}$ and $\{\sigma=\beta\}$.

 \begin{Prop}\label{spectral1}
Let $\sigma=\alpha\mathbbm{1}_{(0,a)}+\beta\mathbbm{1}_{(a,2\pi)}$, where $a\in(0,2\pi)$ and $0<\alpha<\beta<\infty$. Let
\begin{equation}
\nu_1=\min\Bigg\{\frac{\int_0^{2\pi}\sigma|u'|^2\,dt}{\int_0^{2\pi}\sigma u^2\,dt}\,:\,u\in H^1((0,2\pi)),\, u(0)=u(2\pi),\,\int_0^{2\pi}u\,dt=0\Bigg\}.
\end{equation}
Then there exists $\gamma=\gamma\big(\frac{\beta}{\alpha}\big)>\frac{1}{4}$ independent from the parameter $a$, such that $\nu_1>\gamma$.
\end{Prop}
\begin{proof}
The derivative of the functional that defines $\nu_1$ vanishes if and only if
\begin{equation*}
\int_0^{2\pi}\sigma u'v'\,dt=\nu_1\int_0^{2\pi}\sigma uv\,dt,\quad\forall v\in H^1_0((0,2\pi)).
\end{equation*}
We deduce that any optimal  $u$ is a solution of the following system:
\begin{align*}
\begin{cases}
-u''=\nu_1 u&\text{in }(0,a),\\
-u''=\nu_1 u&\text{in }(a,2\pi),\\
u(0)=u(2\pi),\\
\alpha u'(0)=\beta u'(0),\\
\alpha u'(a)=\beta u'(a),\\ 
u \text{ continuous on }(0,2\pi).
\end{cases}
\end{align*}
From the two equations we derive that
\begin{align}
& u(t)=A_1 \cos(\omega t)+A_2\sin(\omega t),\quad \forall t\in(0,a),\\
& u(t)=B_1 \cos(\omega t)+B_2\sin(\omega t),\quad \forall t\in(a,2\pi),
\end{align}
where we have set $\omega=\sqrt{\nu_1}$ and $A_i$, $B_i$ are real constants to be determined. Imposing the continuity conditions in $a$ and $2\pi$, and the transmission conditions in the same points, we get the following system:
\begin{align*}
\begin{cases}
A_1-\cos(2\pi\omega)B_1-\sin(2\pi\omega)B_2=0\\
\cos(a\omega)A_1+\sin(a\omega)A_2-\cos(a\omega)B_1-\sin(a\omega)B_2=0\\
\alpha A_2+\beta\sin(2\pi\omega)B_1-\beta\cos(2\pi\omega)B_2=0\\
-\alpha\sin(a\omega)A_1+\alpha\cos(a\omega)A_2+\beta\sin(a\omega)B_1-\beta\cos(a\omega)B_2=0.
\end{cases}
\end{align*}
Denoting by $A$ the matrix of the coefficients of the previous system, doing some elementary calculations and applying trigonometric identities, we compute
\begin{equation}
\det(A)=\frac{1}{2}(\alpha-\beta)^2[-\cos(2(a-\pi)\omega)+(C+1)\cos(2\pi\omega)-C],
\end{equation}
where 
\begin{equation}
C=C\Big(\frac{\beta}{\alpha}\Big)=\frac{4\frac{\beta}{\alpha}}{\Big(1-\frac{\beta}{\alpha}\Big)^2}>0.
\end{equation}
In order to study $\nu_1$, we need to estimate the first value of $\omega$ that nullifies the following function:
\begin{equation*}
f(\omega,C)=-\cos(2(a-\pi)\omega)+(C+1)\cos(2\pi\omega)-C,\quad\forall\omega>0.
\end{equation*}
We start by finding the   zeros of the function
\begin{equation}
g(\omega)=f(\omega,0)=-2\sin((2\pi-a)\omega)\sin(a\omega),
\end{equation}
which are
\begin{equation*}
\overline{\omega}_k=\frac{k\pi}{a} \quad\text{or}\quad \tilde{\omega}_k=\frac{k\pi}{2\pi-a},\quad\forall k\in\N.
\end{equation*}
Let us assume that $a\in[\pi,2\pi)$ so that $\overline{\omega}_1\leq\tilde{\omega}_1$. It holds that $g\leq 0$ in $\big(0,\frac{\pi}{a}\big)$.

Notice that $\frac{\partial f}{\partial C}=\cos(2\pi \omega)-1$ which is negative. It follows that  $f(\omega,C)\leq g(\omega)$, which yields
\begin{equation*}
f(\omega,C)\leq g(\omega)<0,\quad\forall\omega\in\bigg(0,\frac{\pi}{a}\bigg).
\end{equation*}
Thus by continuity, $f$ must vanish after the first zero of $g$, in other words
\begin{equation*}
\sqrt{\nu_1}=\min_{f(\omega,C)=0}\omega\geq\min_{g(\omega)=0}\omega=\frac{\pi}{a}.
\end{equation*}
If $a\in(0,\pi)$, then $\tilde{\omega}_1<\overline{\omega}_1$ and $g\leq 0$ in $\big(0,\frac{\pi}{2\pi-a}\big)$. With the same argument we get that
\begin{equation*}
\sqrt{\nu_1}=\min_{f(\omega,C)=0}\omega\geq\min_{g(\omega)=0}\omega=\frac{\pi}{2\pi-a}.
\end{equation*}
At this point we have proved that
\begin{equation}
\label{eq1}
\nu_1\geq \min\bigg\{\bigg(\frac{\pi}{a}\bigg)^2,\bigg(\frac{\pi}{2\pi-a}\bigg)^2\bigg\}\geq\frac{1}{4}.
\end{equation}
Let us remark that $0<\omega_1\leq 1$. Indeed, assume first that $\cos(2(a-\pi))\neq 1$. We notice that $f(1,C)=1-\cos(2(a-\pi))> 0$ and on the other hand $f(1/2,C)=-\cos(2(a-\pi))-1-2C < 0$ which proves that $\omega_1 \in(0,1)$. Now if $\cos(2(a-\pi)) = 1$, then 
$$f(\omega,C)=-1+(C+1)\cos(2\pi \omega)-C,$$
which vanishes for $\omega \in \N$ thus $\omega_1=1$ in this case. In any case we have proved that $\omega_1\in (0,1]$.

In other words, we have proved that $\omega_1(a)$ stays in a compact subset of $\R$. Notice that thanks to the bound in \eqref{eq1}, we already know that $\omega_1(a)>1/2$ away from the particular values $a=2\pi$ and $a=0$. However, up to subsequences, if $a\rightarrow 0^+$, we know by compactness that $\omega_1\rightarrow\eta$, for some $\eta\in[0,1]$. Passing to the limit as $a\rightarrow 0^+$ in the eigenvalue equation, we get that
\begin{equation*}
\cos(2\pi\eta)=1, 
\end{equation*}
implying that $\eta\in\{0,1\}$. Since $\eta=\lim_{a\rightarrow 0^+}\omega_1\geq \frac{1}{2}$, it follows that $\eta=1$. The same argument can be applied if we let $a\rightarrow 2\pi^-$. Therefore, there exists $\delta=\delta\big(\frac{\beta}{\alpha}\big)>0$ such that
\begin{equation*}
\omega_1>\frac{3}{4},\quad \forall a\in(0,\delta) \text{ or }\forall a\in(2\pi-\delta,2\pi).
\end{equation*}
If $a\in[\delta,2\pi-\delta]$, then by \eqref{eq1} it holds
\begin{equation}
\omega_1\geq \min\bigg\{\frac{\pi}{\delta},\frac{\pi}{2\pi-\delta}\bigg\}=\frac{\pi}{2\pi-\delta}.
\end{equation}
Thus, choosing $\delta<\frac{2}{3}\pi$, we prove the thesis with $\gamma=\sqrt{\frac{\pi}{2\pi-\delta}}$.
\end{proof}

\begin{Cor}\label{regularity2phases} If $(u,E)$ is a minimizer of the optimal design Problem \eqref{optimalDesign} and let $B_r(x_0)\subset\Omega$. If there exists $r_0>0$ such that $\sharp (\partial E\cap \partial B_r(y))\leq 2$ for all $r<r_0$ and all $y \in \partial E\cap B_{r_0}(x_0)$ then $\partial E$ is smooth in $B_{\frac{r_0}{2}}(x_0)$.
\end{Cor}

\begin{proof} Proposition \ref{spectral1} implies  that for all $y\in \partial E \cap B_{r_0}(x)$ and $r\in (0,r_0)$ we have 
$$\int_{B_r(x)} \sigma_E |\nabla u|^2 \,dx \leq C_0r^{1+\varepsilon},$$
for some $\varepsilon>0$, and we conclude by applying  the theory of almost minimizers for the perimeter~\cite{tamanini}.
\end{proof}




\subsection{The components of $A$ have  mutually quantitative  positive distance}

 \label{SectionDist}
 
 In \cite{Larsen1999} Larsen proved that if $(u,E)$ is a minimizer of the functional \eqref{optimalDesign22}, then for any two components $E_1,E_2$ of $E$ one has 
 $${\rm dist}(E_1,E_2)>0.$$
 
 In this section, we use the uniform rectifiability of $\partial E$ established in Section \ref{uniformRect} to improve the result of Larsen \cite{Larsen1999}. More precisely, we show the validity of Theorem \ref{main2}, whose proof relies on the following lemma. It is a quantitative adaptation of \cite[Lemma 3.1]{Larsen1999}, which we derived using the uniform rectifiability of 
$\dd E$. Actually, the following Lemma contains \cite[Lemma 3.1]{Larsen1999} with a slightly different proof, which is more detailed.

 \begin{Lem}
 \label{LemmaRectangle}
     Let $(u,E)$ be be a minimizer of the functional \eqref{optimalDesign22} and $E_1$, $E_2$ be two connected components such that 
     \begin{equation*}
         \mathrm{dist}(E_1,E_2)^2\leq\delta \min\left\{|E_1|,|E_2|\right\},
     \end{equation*}
    where $\delta\in (0,1)$.   Then for any $R\leq \frac{1}{2\sqrt{\delta}}{\rm dist}(E_1,E_2)$ there exist a ball $B_R(c)$ and two connected components $E'_1$ and $E'_2$ of $E_1\cap B_R(c)$ and $E_2\cap B_R(c)$ respectively such that they are contained in a rectangle. More precisely, up to a rotation, it holds that
    $$E_1'\cup E_2' \subset \{(x,y) \; |\; |y-y_c|\leq h\},$$
    with $h:= C_0(\delta^{\frac{1}{4}}R+R^{\frac{3}{2}})$, where the constant $C_0>0$ depends on the Ahlfors-regularity constant. Moreover, for $i=1,2$, each $E_i'$ contains a point $a_i$ satisfying $|a_i|\leq {\rm dist}(E_1',E_2')$, and $E_i'\cap \partial B_R(c)\neq \emptyset$.
 \end{Lem}
 \begin{proof}
Let $a_1\in E_1$ and $a_2\in E_2$ be such that 
$|a_2-a_1|\leq 2{\rm dist}(E_1,E_2)$. Without loss of generality, we may assume that $\frac{a_1+a_2}{2}=0$. Let 
$$R:=\frac{1}{2\sqrt{\delta}}{\rm dist}(E_1,E_2).$$

By the choice of $R$, it holds that
     \begin{equation}
         E_i\setminus B_R\neq\emptyset,\quad\forall i\in\{1,2\}, \label{nonempty}
     \end{equation}
     because $|B_R|= \pi \frac{1}{4\delta} {\rm dist}(E_1,E_2)^2 \leq  \frac{\pi}{4}  \min\{|E_1|,|E_2|\}< |E_i|$, for $i=1,2$.\\
     \indent Since $E_i$ is connected and contains a point outside $B_R(0)$, this point must be connected to $a_i\in E_i$. For the rest of the proof we will still denote by  $E_1$ and $E_2$ the connected components of $E_1\cap B_R(0)$ and $E_2\cap B_R(0)$ containing respectively $a_1$ and $a_2$ (these components will be the $E_1'$ and $E_2'$ of the statement). Note that from \eqref{nonempty} we deduce that 
$$\partial B_R\cap E_i \neq \emptyset,$$
which is the claim at the end of the statement. 

Let $z_i\in E_i$. Since $E_1$ and $E_2$ are two different connected components, it is true that $\dd E_1$ and $\dd E_2$ separate $z_1$ and $z_2$.  Thus, by  Theorem 14.3 page 123 of \cite{newman}, for any $i\in\{1,2\}$ there exists a connected subset $\Gamma_i$ of $\dd E_i$ such that $\Gamma_i$ separates $z_1$ and $z_2$ and $\mathcal{H}^1(\Gamma_i)<+\infty$. Furthermore, $\Gamma_i$ is arcwise connected.\\
\indent Let $\ell:=[z_1,z_2]$ and let $D$ be the diameter of $B_R$ parallel to $\ell$. We set
\begin{equation}
    \eta:=\mathrm{dist}(\ell,D).
\end{equation}


\begin{figure}[h]
    \centering
    \includegraphics[scale=0.7]{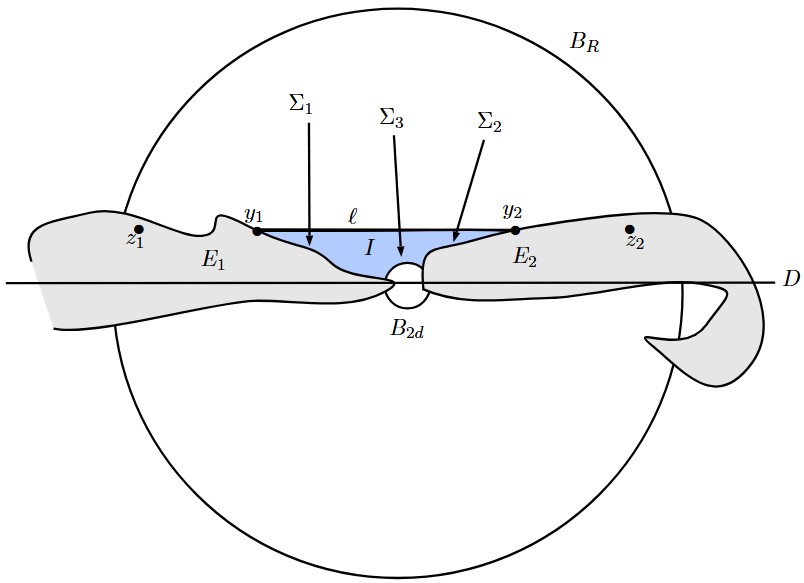}
    \captionsetup{font=small, labelfont={bf, scriptsize}, textfont={it}, width=25cm, margin=2cm}
\caption{\scriptsize The domain noted by $I$ used as a competitor in the proof of Lemma \ref{LemmaRectangle}, following the idea of Larsen \cite{Larsen2003}: if the two connected components $E_1$ and $E_2$ were not flat enough in $B_R$, then we could win a lot of perimeter by adding the domain denoted by $I$ in the picture.}
\end{figure}

Since $\ell\cap\Gamma_i\neq\emptyset$, the following two points exist:
\begin{equation*}
    y_1:=\sup\{t\in[0,1]\,:\, tz_1+(1-t)z_2\in E_1\}\quad\text{and}\quad y_2:=\inf\{t\in[0,1]\,:\, tz_1+(1-t)z_2\in E_2\}.
\end{equation*}
 
Let us denote by $d:=\mathrm{dist}(E_1,E_2)$.  If $\eta\leq 2d $, the conclusion of the lemma holds true. We can assume that $\eta> 2d$. In this case, $B_{2d}\cap\ell=\emptyset$, and $\partial B_{2d}\cap \partial E_i\neq \emptyset$, for $i=1,2$. Thus  there exist three curves $\Sigma_1\subset\dd E_1$, $\Sigma_2\subset\dd E_2$ and $\Sigma_3\subset\dd B_{2d}$ such that the set $\Sigma:=\Sigma_1\cup\Sigma_2\cup\Sigma_3$ is a curve that goes from $y_1$ to $y_2$. By Lemma \ref{touching} below, we know that this curve may possibly have self-intersection points but only by a zero $\mathcal{H}^1$-measure set.\\
 \indent The height bound (Lemma \ref{pythag}) gives
\begin{equation*}
    \mathrm{dist}(x,[y_1,y_2])\leq\sqrt{\frac{\mathcal{H}^1(\Sigma)(\mathcal{H}^1(\Sigma)-|y_2-y_1|)}{2}},\quad\forall x\in\Sigma.
\end{equation*}
Thus, by the triangle inequality and the previous one, we get 
\begin{align}
\label{eqq12}
    \eta=\mathrm{dist}([y_1,y_2],D)
    & \leq \sup_{z\in 
 \Sigma_3}\mathrm{dist}(z,[y_1,y_2])+\sup_{z\in \Sigma_3}\mathrm{dist}(z,D)\leq \sup_{z\in 
 \Sigma}\mathrm{dist}(z,[y_1,y_2])+2d\notag\\
 & \leq \sqrt{\frac{\mathcal{H}^1(\Sigma)(\mathcal{H}^1(\Sigma)-|y_2-y_1|)}{2}}+2d.
\end{align}
At this point we use the $\Lambda$-minimality relation to estimate the right-hand side of \eqref{eqq12}. We denote by $I$ the interior of the Jordan curve $\Sigma\cup[y_1,y_2]$. By the $\Lambda$-minimality of $(u,E)$ with respect to $(u,E\cup I\cup B_{2d})$, we get
\begin{equation*}
    \int_{\Omega}\sigma_E|\D u|^2\,dx + P(E;\Omega)\leq \int_{\Omega}\sigma_{E\cup I\cup B_{2d}}|\D u|^2\,dx + P(E\cup I\cup B_{2d};\Omega) + \Lambda \pi R^2.
\end{equation*}
The inequality can be simplified as
\begin{equation*}
\mathcal{H}^1(\Sigma_1\cup\Sigma_2)\leq \mathcal{H}^1([y_1,y_2])+4\pi d+ \Lambda \pi R^2.
\end{equation*}
By adding $\mathcal{H}^1(\Sigma_3)$ to both sides of the previous inequality, we infer that
\begin{equation*}
\mathcal{H}^1(\Sigma)- \mathcal{H}^1([y_1,y_2])\leq 8\pi d+ \Lambda \pi R^2.
\end{equation*}
Combining \eqref{eqq12} with the previous inequality and using the Ahlfors regularity, we obtain
\begin{align}
    \eta &\leq \sqrt{\mathcal{H}^1(\Sigma)\bigg(4\pi d+\frac{\Lambda}{2}\pi R^2\bigg)}+2d \\
   &\leq C \sqrt{R\bigg(4\pi d+\frac{\Lambda}{2}\pi R^2\bigg)}+2d\\
    &\leq C \sqrt{ \sqrt{\delta} R^2+ R^3}+C\sqrt{\delta} R. \\
    &\leq C(\delta^{\frac{1}{4}}R+R^{\frac{3}{2}}).
\end{align}
The conclusion of the lemma follows by applying Lemma \ref{geometric} below.
\end{proof}

\begin{Lem}\label{geometric} Let $A\subset B:=B_1(0)$ a given set such that $A\cap \partial B \neq \emptyset$ and satisfying the following property: there exists $\eta\leq 1/2$ such that for all $z_1,z_2\subset A$, 
$${\rm dist}([z_1,z_2],D)\leq \eta,$$
where $D$ is the diameter of $B$ parallel to $[z_1,z_2]$. Then, up to a rotation, 
$$
A\subset \{(x,y) \; |\; |y|\leq 3\eta\}.
$$
\end{Lem}

\begin{proof} Let $z_0\in A\cap \partial B$. For any  point $z\in A$ we denote by $D_z$ the the diameter of $B$ parallel to $[z_0,z]$. From our assumption on $A$ we know that for all $z\in A$ it holds,
$${\rm dist}([z_0,z],D_z)\leq \eta.$$
 In particular, $d(z_0,D_z)\leq \eta$ and since $z_0 \in \partial B$, this implies that the angle between $z_0$ and the direction of $D_z$ is small. More precisely, if ${\bf e}_z$ is a unit vector in the direction of $D_z$, then  $$|z_0-\langle z_0 , {\bf e}_z \rangle{\bf e}_z|\leq \eta$$ and, accordingly,
 \begin{equation}
 |z_0-\langle z_0 , {\bf e}_z \rangle{\bf e}_z|^2=1-\langle z_0 , {\bf e}_z \rangle^2\leq \eta^2. \label{estimM}
  \end{equation}
 Let $\theta$ be the angle between the vectors $z_0$ and ${\bf e}_z$, in such a way that $|\langle z_0 , {\bf e}_z \rangle|=\cos(\theta)$. Then we deduce from \eqref{estimM} that
 $$|\sin(\theta)|\leq \eta.$$
 In other words, all the diameters $D_z$, for $z\in A$, must be contained in an angular sector of aperture at most $2\arcsin(\eta)$ around $z_0$. By assuming that $z_0=(1,0)$,  the first vector of the canonical basis  of $\R^2$, we conclude that, for all $z\in A$, $D_z\subset \{(x,y) \; |\; |y|\leq \arcsin(\eta)\}$ and finally, for all $z\in A$, $z \in  \{(x,y) \; |\; |y|\leq \eta + \arcsin(\eta)\}.$ The proposition follows from the elementary inequality, 
 $\arcsin(\eta)\leq 2\eta$, which valid for all $\eta\leq 1/2$.
\end{proof}

We are now ready to proof the main result of this section and of the paper.

\begin{proof}[Proof of Theorem \ref{main2}] We may assume by contradiction that for every $\delta\in(0,1)$ and $\varepsilon>0$ there exist two connected components of $E$ such that 
\begin{equation*}
\min\{|E_1|,|E_2|\}<\varepsilon \quad\text{and}\quad {\rm dist}(E_1,E_2)^2\leq\delta\min(|E_1|,|E_2|).
\end{equation*}
Let $R:= \frac{1}{2\sqrt{\delta}}{\rm dist}(E_1,E_2)$. By Lemma \ref{LemmaRectangle}, there exist a ball $B_R(c)$ and two connected components $E'_1$ and $E'_2$ of $E_1\cap B_R(c)$ and $E_2\cap B_R(c)$  such that  they are contained in a rectangle. More precisely, up to a rotation,
    \begin{equation*}
        E_1'\cup E_2' \subset \{(x,y) \; |\; |y-y_c|\leq h\},
    \end{equation*}
with $h:= C(\delta^{\frac{1}{4}}R+R^{\frac{3}{2}})$, with $C>0$ depending on the Ahlfors regularity constant. Without loss of generality, we may assume that $c=0$.\\
\indent Next, we want to find a curve $\Gamma\subset \partial E'_1 \cap B_R(c)$ which is lower Ahlfors-regular. The easiest way to do so is to find a connected subset of $\partial E_1'$ whose diameter is comparable to the diameter of $E_1'$. For that purpose, we use the following general topological result, which follows directly from the main theorem of \cite{topology1} (see also \cite[Theorem 14.2 p. 123]{newman}): if $U\subset \R^N$ is an open set such that $\R^N\setminus U$ is connected, then $\partial U$ is connected. To simplify the notation, we denote for a moment by $E$ the connected component $E'_1$ . In particular we know that $E\cap B_R(0)\subset \{(x,y) \; | \; |x_2|\leq h\}$.  Let  $U$ be the connected component of $\R^2\setminus E$  containing the point $(0,2h)$. Then $\R^2\setminus U=E\cup_{i\in I} F_i$, where $\{F_i\}_{i\in I}$ is the collection of all the other connected components of $\R^2\setminus E$ different from $U$. In particular, $F_i$ is a closed set, for any $i\in I$. Applying Lemma \ref{topologie} we deduce that $\R^2\setminus U$ is connected. Therefore, from \cite{topology1} we know that $\partial U$ is connected. Furthermore, it is easy to see that $\partial U\subset \partial E$. Let us further show that 
\begin{equation}
{\rm diam}(\partial U) \geq R/4. \label{estimDiam}
\end{equation}
Indeed, from Lemma \ref{LemmaRectangle} we know that $E$ contains a point $a$ satisfying $|a|\leq {\rm dist}(E_1',E_2')=2\sqrt{\delta}R$, and $E\cap \partial B_R\neq \emptyset$. Since $E$ is open and connected, it is pathwise connected. There exists a curve inside $E$ from the point $a$ to a point on $\partial B_R$. This curves stays inside the rectangle $\{(x,y) \; | \; |x_2|\leq h\}$. By consequence, each vertical line $L_t:=t{\bf e}_1 + \R {\bf e}_2$ must intersect $\partial U$ for all $t \in (2\sqrt{\delta}R, R/2)$, and \eqref{estimDiam} follows. We then define $\Gamma=\partial U$ (see Figure 2.). 


\begin{figure}[h]
\centering
 \includegraphics[scale=0.7]{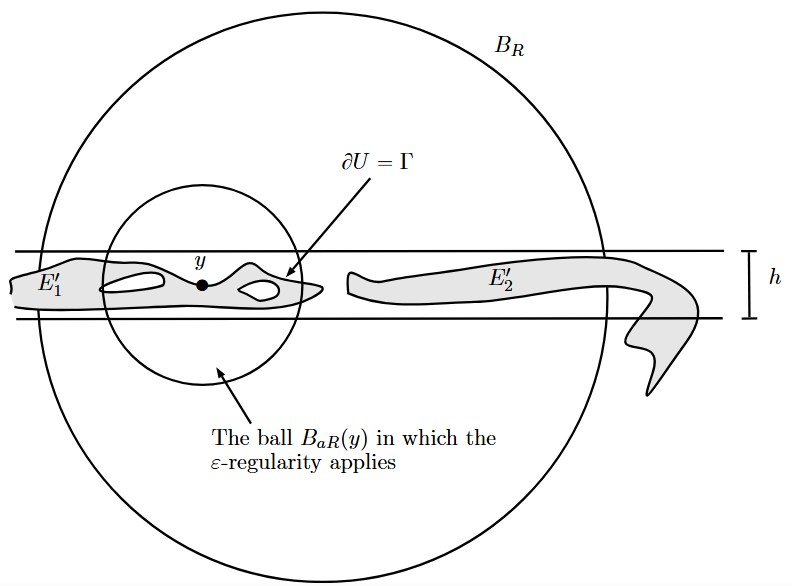}
 \captionsetup{font=small, labelfont={bf, scriptsize}, textfont={it}, width=25cm, margin=2cm}
\caption{\scriptsize In the proof of Theorem \ref{main2}, we get a contradiction by finding a ball in which $\partial E$ should be a smooth surface thanks to the  $\varepsilon$-regularity result,  which prevents a component to be shrinked in a thin rectangle.}
\end{figure}

 
Since    $\Gamma$ is connected, it  satisfies the following lower Ahlfors regularity property:
$$\mathcal{H}^{1}(\Gamma \cap B_r(x))\geq r, \quad\forall x \in \Gamma, \,\forall r \in\bigg(0,\frac{R}{4}\bigg).$$

Then by Lemma \ref{porosity} there exist two constants $a\in(0,1)$, $\varepsilon_0>0$ and a ball $B_{a R}(y)\subset B_R(0)$ such that 
\begin{equation*}
    \mathbf{e}(y,a R)+\omega(y,aR)+a R\leq\varepsilon_0,
\end{equation*}
for $R$ sufficiently small. By $\varepsilon$-regularity we have that $\dd E'_1\cap B_{\frac{a R}{2}}(y)$ is a $C^{1,\gamma}$-hypersurface. We can choose $\delta<\big(\frac{a}{4C}\big)^4$ and $\varepsilon<\big(\frac{a}{4C}\big)^4$ so that the radius of the ball $B_{\frac{a R}{2}}(y)$ is greater than the height $h$ of the rectangle, that is
$C\big(\delta^{\frac{1}{4}}R+R^{\frac{3}{2}}\big)<\frac{a R}{2}$. Indeed, 
\begin{equation*}
    C\Big(\delta^{\frac{1}{4}} R+R^{\frac{3}{2}}\Big)< \frac{C_0R}{4}+C\frac{(\min\{|E_1|,|E_2|\})^{\frac{1}{4}}}{4}<\frac{aR}{4}+\frac{C\varepsilon^{\frac{1}{4}}}{4}=\frac{aR}{2},
\end{equation*}
which concludes the proof. Indeed, this is clearly a  contradiction with the fact that   $E_1' \cap B_R$ was supposed to be totally contained in the rectangle.
\end{proof}



The following lemma has been used in the proof of Lemma \ref{LemmaRectangle} and it holds when $N=2$. It states that, under some mild regularity assumptions on a $BV$ set $\Omega$, the boundaries of connected components cannot touch by a positive $\mathcal{H}^{1}$-measure set. 

\begin{Lem}\label{touching}
Let $\Omega\subset \R^2$ be a set of finite perimeter satisfying a lower Ahlfors-regularity inequality and such that 
\begin{equation}
\label{assumm}
\mathcal{H}^{1}(\partial \Omega \setminus \partial^*\Omega)=0.
\end{equation}
 Let $A$ and $A'$ be two connected components of $\Omega$ and let $x_0\in \partial A \cap \partial A'$. Then $x_0$ does not belong to $\partial^*\Omega$. As a consequence,
$$\mathcal{H}^{1}(\partial A \cap \partial A')=0.$$
\end{Lem}

\begin{proof} Assume by contradiction that  $x_0 \in \partial^*\Omega \cap \partial A \cap \partial A'$. We know that $\partial^*\Omega$ admits an approximative tangent plane $P_0$ at point  $x_0$. From \eqref{assumm} we deduce that  $P_0$ is actually an approximative tangent plane for   $\partial \Omega$ as well. Then since $\partial \Omega$ satisfies a lower Ahlfors-regularity inequality, it is classical to see that this approximative tangent plane is actually a true tangent plane. Let us write more details about this last fact.

We assume without loss of generality that $P_0=\{x_2=0\}$ and let $\varepsilon>0$ be fixed. Since $P_0$ is an approximative tangent plane we know that  
\begin{equation}
\lim_{r\to 0} \frac{1}{r}\mathcal{H}^1(\partial \Omega \cap B_r(x_0)\cap \{|x_2|\geq  r\varepsilon \} ) \to 0.\label{avance}
\end{equation}
Now we claim that there exists $r_0>0$ such that,  for all $r\leq r_0$, 
$$\partial \Omega \cap B_r(x_0)\cap \{|x_2|>  r\varepsilon \} = \emptyset.$$
Otherwise, if $z\in \partial \Omega \cap B_r(x_0)\cap \{|x_2|> r\varepsilon \}$, since $\partial \Omega$ is lower Ahlfors regular then 
$$\mathcal{H}^{1}(\partial \Omega \cap B_{\varepsilon r}(z))\geq C_A \varepsilon r,$$
Which would easily contradict \eqref{avance} for $r$ small enough.\\
\indent We may assume that for $\varepsilon>0$ we can find $r_0>0$ such that for any $r\in(0,r_0)$ we have 
\begin{equation}
\label{eqq14}
    \dd\Omega\cap B_r(x_0)\subset T_r(\varepsilon):=\{|x_2|\leq \varepsilon|x_1|\}.
\end{equation}
Without loss of generality we may assume that $x_0=0$. We distinguish two cases.\\
\indent \textbf{Case 1:} there exists $r\in(0,r_0)$ such that 
\begin{equation*}
    (B_r\setminus T_r(\varepsilon))\cap A\neq\emptyset\quad\text{and}\quad (B_r\setminus T_r(\varepsilon))\cap A'\neq\emptyset.
\end{equation*}
Let $z\in (B_r\setminus T_r(\varepsilon))\cap A$ and $z'\in (B_r \setminus T_r(\varepsilon))\cap A'$. 
Let us define
\begin{equation*}
    B_r\setminus T_r(\varepsilon)=\big(B_r\setminus T_r(\varepsilon)\cap\{x_2>0\}\big)\cup \big(B_r\setminus T_r(\varepsilon)\cap \{x_2<0\}\big)=: T^+\cup T^-.
\end{equation*}
Notice that if $z\in T^+$, then necessarily $T^+\subset A$. Indeed, if $y\in T^+$ is any other point, then the segment $\overline{yz}$ is contained in $T^+$ because $T^+$ is convex. Since $z$ belongs to the connected component $A$, it holds that $y\in A$, thus proving $T^+\subset A$. The same assertion holds for $A'$ so that  the following two alternatives must hold:
\begin{enumerate}
    \item $T^+\subset(B_r\setminus T_r(\varepsilon))\cap A$ and $T^-\subset(B_r\setminus T_r(\varepsilon))\cap A'$.
    \item $T^+\subset(B_r\setminus T_r(\varepsilon))\cap A'$ and $T^-\subset(B_r\setminus T_r(\varepsilon))\cap A$.
\end{enumerate}

\indent In both cases, it follows that there exists a positive constant $C$ such that
\begin{equation}
    |\Omega\cap B_r|\geq \pi r^2-C\varepsilon r^2,
\end{equation}
which is a contradiction if $\varepsilon<\frac{\pi}{2 C}$, being $x_0=0\in\dd^*\Omega$.\\
\indent \textbf{Case 2:} for any $r\in(0,r_0)$ it holds
\begin{equation*}
    (B_r\setminus T_r(\varepsilon))\cap A=\emptyset\quad\text{or}\quad (B_r\setminus T_r(\varepsilon))\cap A'=\emptyset.
\end{equation*}
We assume without loss of generality that $(B_r\setminus T_r(\varepsilon))\cap A=\emptyset$. We take a point $z\in T_r(\varepsilon)\cap A$. Without loss of generality, we can assume that $s_0:=\pi_1(z)>0$. Since $A$ is connected, and $z,0\in A$, we deduce that $\partial B_{s}\cap A \neq \emptyset$, for all $s\in(0,s_0)$. Let $z_s$ be a point in $\partial B_{s}\cap A$. Since $(B_r\setminus T_r(\varepsilon))\cap A=\emptyset$, we know that $z_s\in  T_r(\varepsilon)$.

Since $0\in\dd A$, for all $s<s_0$ there exists a point $z'_s\in A\cap B_{\varepsilon s}(x_0)$ such that $\pi_1(z'_s)\leq \pi_1(z_s)$. Let $\gamma$ be a curve connecting $z_s$ and $z'_s$ in $A$, which exists because $A$ is a connected open set, thus arc-wise connected.\\
\indent  We define the vertical line
\begin{equation*}
    L_t:=\{(t,y)\,:\,y\in\R\},\quad\text{for }t\in\R.
\end{equation*}

Since $\dd A\cap B_s(x_0)\subset\dd \Omega\cap B_r(x_0)\subset T_r(\varepsilon)$ and $L_t$ meets interior points of $A$ for $t\in[\pi_1(z'_s),\pi_1(z_s)]$, we have that
\begin{equation*}
    \sharp(\dd A\cap L_t)\geq 2,\quad\forall t\in[\pi_1(z'_s),\pi_1(z_s)].
\end{equation*}
Accordingly, by the coarea formula we get
\begin{align*}
    \mathcal{H}^1(\dd A\cap B_s(x_0))&\geq \int_{\dd A\cap B_s(x_0)}\sqrt{1-\prodscal{\nu_{A}}{e_n}^2}\,d\mathcal{H}^1 \\&\geq\int_{s\varepsilon}^{s/\sqrt{1+\varepsilon^2}}\mathcal{H}^0(\dd A\cap\{\pi_1=t\})\,dt\geq 2s\bigg(\frac{1}{\sqrt{1+\varepsilon^2}}-\varepsilon\bigg)  \\
    &\geq 2s(1-2\varepsilon),
\end{align*}
for any $s\in(0,s_0)$. Moreover we know that for all $s\in (0,s_0)$,
$$\partial \Omega \cap \partial B_s(x_0)\cap \{\pi_1(x)<0\} \neq \emptyset,$$ 
so that, all together,
\begin{equation*}
    \mathcal{H}^1(\dd\Omega\cap B_s(x_0))\geq s+2s(1-2\varepsilon)=s(3-4\varepsilon),
\end{equation*}
which implies the following contradiction:
\begin{align*}
    \liminf_{s\rightarrow0}\frac{\mathcal{H}^1(\dd^*\Omega\cap B_s(x_0))}{2s}=\liminf_{s\rightarrow0}\frac{\mathcal{H}^1(\dd\Omega\cap B_s(x_0))}{2s}\geq \frac{3}{2}-2\varepsilon>1,
\end{align*}
for $\varepsilon<\frac{1}{4}$.

\end{proof}

\begin{Rem} The Lemma is false in higher dimensions: consider the $1D$-curve $\Gamma \subset \R^N$ defined by 
$$\Gamma:=\{x_1=x_2=\dots =x_{N-2}=0 \text{ and } x_N=x_{N-1}^2\}.$$
Then consider $G$ as being a very small open neighborhood of $\Gamma$, with the property that $G\subset \{x_N>0\}$ and
$$\lim_{r\rightarrow 0^+}\frac{|G\cap B_r(0)|}{r^N}= 0.$$
Then we define $\Omega:= P^+ \cup G$, where $P^+=\{x_N<0\}$. The two connected components of $\Omega$ are $P^+$ and $G$, for which their boundaries meet at the origin, which is a point that belongs to $\partial^*\Omega$.  However, we don't know if the assertion $\mathcal{H}^{N-1}(\partial A \cap \partial A')=0$ still holds in higher dimensions for two connected components of a set of finite perimeter $\Omega$ satisfying $\mathcal{H}^{N-1}(\partial \Omega \setminus \partial^*\Omega)=0$, and such that $\partial \Omega$ satisfies a lower Ahlfors-regularity inequality. 
\end{Rem}


\section{Appendix}

The following lemma is a standard  height bound which is taken from   \cite[Lemma 6.3]{LemBaba} for the case of injective curves. Here  we adapted the argument straightforwardly for a curve which is possibly non injective anymore, but for which the set of self-intersection points has zero measure.  

 \begin{Lem}\label{pythag} Let $\gamma:[0,1]\to \R^2$ be a curve with endpoints $z=\gamma(0)$ and $z'=\gamma(1)$, with image $\Gamma:=\gamma([0,1])$. We assume that $\gamma$ is almost injective in the sense that, defining $$Z:=\{t\in [0,1] \,:\,\exists t'\neq t \; \text{ s.t. }\; \gamma(t)=\gamma(t')\},$$
it holds that $\mathcal{H}^1(\gamma(Z))=0$.
It follows that 
  \begin{equation}
 {\rm dist}(y,[z,z'])^2 \leq \frac{\mathcal{H}^1(\Gamma) \big(\mathcal{H}^1(\Gamma)-|z'-z|\big)}{2},  \quad\forall y\in \Gamma. \label{estimgeom}
 \end{equation} 
 \end{Lem}
 
\begin{proof}
Let $\bar y$ be a maximizer of the function $y \in \Gamma \mapsto {\rm dist}(y,[z,z'])$, i.e., $\bar y$ is the most distant point in $\Gamma$ from the segment $[z,z']$, and define $d=:{\rm dist}(\bar y,[z,z'])$. Let us consider the point $y' \in \R^2$ making $(z,z',y')$ an isosceles triangle with height $d$. Denoting by $a:=|z-z'|/2 $ and $L := |y'-z|$, according to Pythagoras Theorem, we have 
 $$d^2 = L^2-a^2=(L -a) (L +a).$$

And since  by assumption $\gamma$ is almost injective (i.e. $\mathcal{H}^1(Z)=0$), then $\mathcal{H}^1(\Gamma)\geq |z-\bar y| + |\bar y - z'|\geq 2L$ and $\mathcal H^1(\Gamma)\geq |z-z'|$ so that
\[d^2\leq \frac{1}{4}\left(\mathcal{H}^1(\Gamma )-|z-z'|\right) \left(\mathcal{H}^1(\Gamma )+|z-z'|\right)\leq  \frac{\mathcal{H}^1(\Gamma )\big(\mathcal{H}^1(\Gamma )-|z-z'|\big)}{2} ,\]
which proves \eqref{estimgeom}.
 \end{proof}

\begin{Lem}\label{topologie} Let $A\subset \R^N$ be an open and connected set, and let $\{F_i\}_{i\in I}$ be a family of closed connected sets such that $F_i \cap \partial A \neq\emptyset$ for all $i \in I$. Then the set
$$A\cup \bigcup_{i \in I} F_i$$
is connected.
 \end{Lem}
\begin{proof}Let us denote by
$$E:=A\cup \bigcup_{i \in I} F_i,$$
and assume  that 
$$E \subset U_1 \cup U_2$$
where $U_1$ and $U_2$ are two disjoint open and connected subsets of $\R^N$. To prove that $E$ is connected, it is enough to prove that $E\subset U_1$ or $E\subset U_2$. 
Let $i \in I$ be fixed for a moment. Then, by assumption, $F_i \subset U_1 \cup U_2$ but since $F_i$ is connected we deduce that $F_i\subset U_1$ or $F_i\subset U_2$. Let us assume that $F_i\subset U_1$. Let $x_0 \in F_i\cap \partial A$, which is assumed to be non empty. Since $U_1$ is open and $x_0\in \partial A$, we actually infer that $U_1 \cap A\neq \emptyset$. But since $A\subset U_1\cup U_2$, and since $A$ is connected, we conclude that $A\subset U_1$. In other words we have proved that $A\cup F_i \subset U_1$.  Now for any other $j \in I$, arguing similarly we deduce that $A\cup F_j$ is also contained  in one of the two open sets, that actually must be the same $U_1$ because $A$ is already known to be contained in $U_1$. All in all we have proved that $E\subset U_1$, as desired, and this achieves the proof.
 \end{proof}

\bibliographystyle{plain}
\bibliography{biblio_OD}

\end{document}